\documentclass[11pt]{article}

\usepackage[margin=1in]{geometry}
\usepackage{amsmath,amssymb,mathtools,amsthm}
\usepackage{booktabs,array}
\usepackage{microtype}
\usepackage{enumitem}
\usepackage{xcolor}
\usepackage[normalem]{ulem}
\usepackage{hyperref}

\hypersetup{colorlinks=true,linkcolor=black,
  citecolor=black,urlcolor=black}
\setlist{nosep,leftmargin=*}
\newtheorem{theorem}{Theorem}
\newtheorem{lemma}{Lemma}
\newtheorem{proposition}{Proposition}

\newcommand{\inner}[2]{\left\langle #1,#2\right\rangle}
\newcommand{\proj}{\operatorname{proj}}
\newcommand{\diag}{\operatorname{diag}}
\newcommand{\one}{\mathbf 1}
\newcommand{\R}{\mathbb R}

\title{Box-Cone Representations of Bound-Factor RLT Constraints in Semidefinite Relaxations}
\author{Katherine Charbonneau$^{1}$ \qquad Hao Hu$^{1}$ \qquad
Renata Sotirov$^{2}$\\[0.5em]
\small $^{1}$School of Mathematical and Statistical Sciences,
Clemson University, Clemson, South Carolina, USA\\
\small $^{2}$Department of Econometrics and Operations Research,
Tilburg University, Tilburg, The Netherlands\\[0.25em]
\small \texttt{kcharbo@clemson.edu}, \texttt{hhu2@clemson.edu},
\texttt{r.sotirov@uvt.nl}}
\date{\today}

\begin{document}
\maketitle

\begin{abstract}
Linear inequalities can generate large constraint systems in semidefinite
relaxations of binary quadratic programs, making the quadratic subproblems
of splitting methods expensive to solve.  We develop a box-cone representation
of bound-factor reformulation-linearization technique (RLT) constraints
in semidefinite relaxations of binary quadratic programs.  These constraints
arise by multiplying each original linear inequality by every binary variable
and its complement and linearizing the resulting products.  Our key observation is that all
products generated by one inequality share a common slack and admit an exact
representation in the standard cone over the unit box.  This box-cone
representation processes the inequalities through independent explicit
projections, without introducing a separate nonnegative slack for each
scalar RLT inequality.  We illustrate the computational use of the representation within standard
two-block alternating direction method of multipliers (ADMM).
The resulting subproblems require linear solves with a reusable
factorization and explicit cone projections, without enlarging the
semidefinite matrix.  Experiments on 257 binary quadratic programs compare
this implementation with SDPNAL+.  A second application groups Max-Cut
triangle inequalities into box cones and uses the resulting bounds within
an existing branch-and-bound framework.  Tests on 220 graphs show
competitive performance with MADAM using only triangle cuts.  Together,
these applications illustrate the representation's usefulness as a
building block for first-order methods.
\end{abstract}


\section{Introduction}
\label{sec:introduction}

We develop a box-cone model of the bound-factor reformulation-linearization
technique (RLT) constraints in semidefinite programming (SDP) relaxations of
binary quadratic programs.  The model represents all bound-factor products
generated by one inequality through a single box cone, preserving the
relaxation while exposing a compact coupling and an explicit projection.
Consider the linearly constrained binary quadratic program
\begin{equation}
\label{eq:bqp}
  \min\{x^\top Qx+q^\top x:Ax\le b,\ Gx=h,\ x\in\{0,1\}^n\},
\end{equation}
where $Q\in\mathbb S^n$, $q\in\R^n$,
$A\in\R^{m\times n}$, $b\in\R^m$, $G\in\R^{p\times n}$, and
$h\in\R^p$; the data are unrestricted in sign, and $Q$ need not be
positive semidefinite.  Either $m$ or $p$ may be zero.  Applications of
\eqref{eq:bqp} include quadratic assignment \cite{oliveira-wolkowicz-xu},
graph partitioning \cite{demeijer-sotirov-wiegele-zhao}, and quadratic
knapsack problems \cite{pisinger-qkp}.

SDP relaxations provide bounds for these
problems and can support exact solution methods such as branch-and-bound.
They represent the variables and their pairwise products in a lifted matrix
of order $n+1$, impose positive semidefiniteness and the binary identities,
and may impose entrywise nonnegativity.  Bound-factor RLT \cite{sherali-adams} strengthens an SDP or linear
programming relaxation by multiplying each linear constraint by the binary
variables and their complements and linearizing the resulting products.  This can substantially strengthen the
bound, but also creates a large system of linear constraints: each original
inequality generates $2n$ bound-factor RLT inequalities in addition to the
original inequality.  Processing these constraints can become a major cost of
solving the relaxation.

In the alternating direction method of multipliers (ADMM), the treatment
of these inequalities determines the cost of the subproblems.  Existing
approaches include iterative polyhedral projection and scalar
slack-variable reformulations.  Retaining coupled inequalities in a single
polyhedral block gives a constrained quadratic projection subproblem;
for example, de Meijer et al.~\cite{demeijer-sotirov-wiegele-zhao,demeijer-siebenhofer-sotirov-wiegele}
use Dykstra's algorithm for projections involving additional cuts.
Cutting-plane schemes can also limit the inequalities imposed at each stage
\cite{demeijer-siebenhofer-sotirov-wiegele}, while related methods exploit
disjoint constraint supports to organize independent projections \cite{demeijer-sotirov-qccp,sinjorgo-sotirov-vera}. 
Slack-variable approaches can expose simple bound projections
\cite{wiegele-zhao,battista-desantis}, but retain auxiliary coordinates
for the individual scalar inequalities.  Such formulations can exploit
additional structure: Battista and De Santis \cite{battista-desantis}
avoid explicitly forming the enlarged slack-based SDP.

The box-cone model exploits the common slack shared by the bound-factor
products generated by each inequality.  For a row $a_r^\top x\le b_r$,
let $t=b_r-a_r^\top x$ and let $u_j$ denote the linearized product $tx_j$.
The product $t(1-x_j)$ becomes $t-u_j$, so the original inequality and its
bound-factor RLT constraints are exactly $t\ge0$ and $0\le u_j\le t$
for $j=1,\ldots,n$.  These conditions define the standard box cone
\cite{odonoghue-box}.  The inequalities in $Ax\le b$, together with their bound-factor RLT
constraints, admit a Cartesian-product representation with one box cone
and an independent projection per row.

This representation combines the two product families in one matrix image
and exposes a separable projection computed by sorting followed by
coordinatewise projection onto an interval.  The image has $m(n+1)$ coordinates,
compared with $m(2n+1)$ scalar slack coordinates in a formulation that
assigns one to every original and product inequality, reducing
storage for the auxiliary variables and their associated multipliers.

The RLT families and the box cone are established objects
\cite{sherali-adams,odonoghue-box}. Matrix formulations of the bound-factor inequalities are also
established~\cite{dey-kocuk}. The box-cone projection is available in
SCS~\cite{scs}, as documented in \cite[Section~6.2]{odonoghue-box}.  Our contribution is to identify and exploit
their shared-slack structure in a formulation that combines the compact
cone image with an explicit affine update.  The representation exposes independent box-cone projections and a
structured linear coupling, allowing the bound-factor inequalities
to be processed together without an iterative polyhedral projection.
These properties make it a useful building block for designing
efficient first-order methods for SDP--RLT relaxations.

We demonstrate this construction within standard two-block ADMM
\cite{boyd-admm,wen-goldfarb-yin}.  The splitting separates the box-cone
image from the normalization and arrow equations and incorporates the
usual facial reduction for equality-generated RLT.  Each iteration then
uses explicit cone projections and the reusable affine solve, without
an iterative polyhedral projection or a separate slack coordinate for
every scalar RLT inequality.  The explicit inequality projections act
in the cone-image space, rather than on the full RLT feasible set in the
lifted matrix space.  The role of formulation is analogous to that of
facial reduction in the QAP method of Oliveira, Wolkowicz, and Xu
\cite{oliveira-wolkowicz-xu}, where it exposes inexpensive semidefinite
and polyhedral updates.  Related ADMM and Peaceman--Rachford
applications also exploit equality-induced facial reduction
\cite{hu-sotirov-qspp,demeijer-sotirov-qccp,hu-sotirov-wolkowicz,xu-hu-qvdp,graham-hu-im-li-wolkowicz,li-pong-sun-wolkowicz,burkowski-im-wolkowicz}.  Numerical experiments compare our ADMM
implementation with SDPNAL+ \cite{yang-sun-toh} on the same SDP--RLT formulation on each instance. 
A second application uses the shared-bound structure of Max-Cut triangle
inequalities to construct box-cone bounds within the branch-and-bound
framework of MADAM \cite{hrga-povh-madam}.  This tests the modeling idea
with dynamically selected cuts as well as with a fixed RLT system.

\paragraph{Organization.}

Section~\ref{sec:model} defines the relaxation precisely.
Section~\ref{sec:star} develops the box-cone representation, its projection,
and its storage implications.  Section~\ref{sec:admm} applies standard
two-block ADMM, including the affine projection and the usual facial
reduction for equality-generated RLT.  Section~\ref{sec:stronger-rlt}
discusses stronger RLT systems.  Implementation details for the fixed SDP--RLT relaxation are given in Section~S1 of the supplement.  Section~\ref{sec:experiments} presents
the SDPNAL+ comparison and the Max-Cut application within MADAM.  Section~\ref{sec:conclusions} discusses the scope of the
contribution and future work on low-rank methods.  The appendix gives the
box-cone representation used in the Max-Cut application.

\section{Model and scope}
\label{sec:model}

We construct an SDP--RLT relaxation of problem~\eqref{eq:bqp} that combines
positive semidefiniteness, entrywise nonnegativity, and RLT constraints to
obtain strong lower bounds.  Bound-factor RLT constraints can strengthen these SDP bounds, but their
large number can make the relaxation expensive to solve.  Our key result enables their efficient treatment through an exact
rowwise representation with an explicit projection.  We now formulate the
relaxation in the notation used for this representation.

For a positive integer $d$, $\mathbb S^d$ denotes the real symmetric
matrices of order $d$; $\mathbb S_+^d$ and $\mathbb S_{++}^d$ denote their
positive semidefinite and positive definite subsets, respectively. Introduce the Shor matrix $Y\in\mathbb S^{n+1}$ and objective
matrix $C\in\mathbb S^{n+1}$ by
\[
Y=\begin{bmatrix}1&x^\top\\x&X\end{bmatrix},\qquad
C=\begin{bmatrix}0&q^\top/2\\q/2&Q\end{bmatrix}.
\]
Here $X\in\mathbb S^n$ represents the pairwise products $x_ix_j$, so
the lifted objective is $\inner{C}{Y}$. We index the rows and columns of
$Y$ by $0,\ldots,n$. In the exact binary lifting,
$Y=[1;x][1;x]^\top$ has rank one. The SDP relaxation drops this rank-one
requirement and imposes $Y\succeq0$, entrywise
nonnegativity $Y\ge0$, the normalization $Y_{00}=1$, and the arrow equations
$Y_{jj}=Y_{0j}$ encoding the binary identities $x_j^2=x_j$ for $j=1,\ldots,n$.

The RLT constraints are obtained by forming valid products of each linear constraint
with $x_j$ and $1-x_j$ and then replacing every product $x_ix_j$ by $X_{ij}$.
Let $a_r^\top$ be row $r$ of $A$ and $g_\ell^\top$ row $\ell$ of $G$.
Introduce the two augmented row operators
\[
 H=[b,-A]\in\R^{m\times(n+1)},\qquad
 F=[h,-G]\in\R^{p\times(n+1)}.
\]
Left multiplication of $Y$ by $H$ collects all inequality slacks and their lifted
products in one matrix:
\[
 HY=\begin{bmatrix}b-Ax & bx^\top-AX\end{bmatrix}.
\]
The zeroth column is the slack vector, and column $j$ contains the lifted products with
$x_j$. For row~$r$, denote the slack and lifted products by
\[
 \sigma_r(Y):=(HY)_{r0}=b_r-a_r^\top x,
 \qquad
 \pi_{rj}(Y):=(HY)_{rj}=b_rx_j-a_r^\top X_{\cdot j}.
\]
For the inequality $a_r^\top x\le b_r$, the valid products are
$(b_r-a_r^\top x)x_j\ge0$ and
$(b_r-a_r^\top x)(1-x_j)\ge0$.  Linearization gives
$\pi_{rj}(Y)\ge0$ and $\sigma_r(Y)-\pi_{rj}(Y)\ge0$, respectively.
Thus the original
inequality and its bound-factor RLT constraints are exactly
\begin{equation}
\label{eq:inequality-rlt}
 \sigma_r(Y)\ge0,\qquad
 0\le \pi_{rj}(Y)\le \sigma_r(Y)
 \quad(r=1,\ldots,m,\ j=1,\ldots,n).
\end{equation}

For an equality $g_\ell^\top x=h_\ell$, linearizing
$(h_\ell-g_\ell^\top x)x_j=0$ gives
$h_\ell x_j-g_\ell^\top X_{\cdot j}=0$.
Together with the original equalities, these constraints are represented by
\begin{equation}
\label{eq:equality-rlt}
 FY=0.
\end{equation}
The zeroth column gives $h-Gx=0$, and column $j$ gives the linearized
products with $x_j$.
The products with $1-x_j$ are redundant: each is the difference between the
zeroth-column equation and the corresponding column-$j$ equation.

Combining these constraints gives the doubly nonnegative (DNN)--RLT relaxation:
\begin{equation}
\label{eq:mixed-rlt}
\begin{aligned}
 \min_Y\quad&\inner{C}{Y}\\
 \text{s.t.}\quad
 &Y_{00}=1,\quad Y_{jj}=Y_{0j}\quad(j=1,\ldots,n),\\
 &Y\ge0,\quad Y\succeq0,\\
 &FY=0,\\
 &\sigma_r(Y)\ge0,\quad 0\le \pi_{rj}(Y)\le \sigma_r(Y)
       \quad(r=1,\ldots,m,\ j=1,\ldots,n).
\end{aligned}
\end{equation}
Here $Y\ge0$ is entrywise. The DNN--RLT relaxation~\eqref{eq:mixed-rlt} contains the original
linear constraints, equality-generated RLT, and all bound-factor RLT constraints in
\eqref{eq:inequality-rlt}.  It does not contain off-diagonal McCormick
inequalities or products between pairs of inequality slacks.  Subsection~\ref{sec:mccormick-row-star}
describes an optional full-McCormick extension, which is not used in the
SDPNAL+ comparison.  The separate Max-Cut application in
Subsection~\ref{sec:maxcut-experiments} uses selected triangle inequalities
in correlation-matrix coordinates.  Subsection~\ref{sec:inequality-products} describes the
corresponding representation of products between inequality slacks and
explains why it is less attractive computationally.

\section{Box-cone representation and projection}
\label{sec:star}

\subsection{The box-cone representation and box-cone projection}

Define the polyhedral cone
\[
 \mathcal S
 :=\operatorname{cone}\bigl(\{1\}\times[0,1]^n\bigr)
 =\{(t,u)\in\R\times\R^n:t\ge0,\ 0\le u\le t\one\}.
\]
This is the standard box cone \cite{odonoghue-box}, or equivalently the homogenization of the unit cube
\cite{burer-dong}.
Here~$t$ is a common nonnegative upper bound on the coordinates of $u$.
If the model has $m$ inequalities, define the product of their box cones by
\[
 \mathcal S_m
 :=\underbrace{\mathcal S\times\cdots\times\mathcal S}_{m\text{ factors}}
 \subseteq\R^{m\times(n+1)},
\]
with one matrix row identified with each factor.

\begin{theorem}[Box-cone identity]
\label{thm:row-star}
For any lifted matrix $Y$, the bound-factor RLT constraints
\eqref{eq:inequality-rlt} are equivalent to the single matrix-cone constraint
\[
       HY\in\mathcal S_m.
\]
\end{theorem}

\begin{proof}
By definition, $(HY)_{r0}=\sigma_r(Y)$ and
$(HY)_{rj}=\pi_{rj}(Y)$.  Membership of row $r$ in $\mathcal S$ is
therefore precisely $\sigma_r(Y)\ge0$ and
$0\le\pi_{rj}(Y)\le\sigma_r(Y)$ for every $j$.
\end{proof}

The shared upper bound is the essential structure.  Keeping only the
lower-bound-factor products would give the orthant constraint $HY\ge0$.
Including the upper-bound-factor products adds
$(HY)_{rj}\le(HY)_{r0}$, so all $2n+1$ scalar
inequalities generated by one original inequality become one cone constraint with
an explicit Euclidean projection.

\begin{samepage}
\begin{proposition}[Closed-form projection onto the box cone]
\label{prop:star-projection}
Given $c=(c_0,c_1,\ldots,c_n)\in\R^{n+1}$, its Euclidean projection
$(t^\star,u^\star)=\proj_{\mathcal S}(c)$ solves
\[
 \min_{\substack{t\ge0,\\0\le u\le t\one}}
 \frac12(t-c_0)^2+\frac12\sum_{j=1}^n(u_j-c_j)^2.
\]
Define
\[
 \phi(t)=t-c_0-\sum_{j=1}^n\max\{c_j-t,0\},\qquad t\ge0.
\]
The projection is obtained by first setting
\[
 t^\star=
 \begin{cases}
  0,&\phi(0)\ge0,\\
  \text{the unique positive root of $\phi(t)=0$},&\phi(0)<0,
 \end{cases}
\]
and the remaining coordinates are
\[
 u_j^\star=\min\{\max\{c_j,0\},t^\star\},
 \qquad j=1,\ldots,n.
\]
\end{proposition}
\end{samepage}

\begin{proof}
For fixed $t$, scalar projection gives
$u_j(t)=\min\{\max\{c_j,0\},t\}$.
Because $t\ge0$, this implies
\[
 (u_j(t)-c_j)^2
 =\max\{-c_j,0\}^2+\max\{c_j-t,0\}^2.
\]
Substituting these coordinates reduces the problem to
\[
 \min_{t\ge0}\ \psi(t),\qquad
 \psi(t)=\frac12(t-c_0)^2
 +\frac12\sum_{j=1}^n\max\{c_j-t,0\}^2
 +\frac12\sum_{j=1}^n\max\{-c_j,0\}^2.
\]
The last sum is constant, and $\psi$ is a strictly convex
piecewise-quadratic function.  For $t>0$, its derivative is
$\psi'(t)=\phi(t)$, while its right derivative at zero is $\phi(0)$.
Moreover, $\phi$ is continuous, strictly increasing, and tends to infinity.
Therefore the minimizer is zero when $\phi(0)\ge0$ and otherwise equals the
unique positive root of $\phi$.
Substituting this value into the interval projection formula gives the stated $u^\star$.
\end{proof}

The scalar $t^\star$ can be computed by sorting the positive $c_j$ and testing
the associated prefix averages.  Because $\mathcal S_m$ is a Cartesian
product of $m$ copies of $\mathcal S$, projection onto $\mathcal S_m$ applies
the single-row formula independently to each row.  The resulting cost is
$O(mn\log n)$ comparisons and $O(mn)$ arithmetic operations, and the row
projections can be evaluated together.

\subsection{Storage and operator structure}

The cone image $HY$ stores each original slack once and each lower-bound
product once.  The upper-bound products are obtained by subtracting the
latter from the former.  Thus one image or associated multiplier has
$m(n+1)$ entries, whereas a separate scalar coordinate for every inequality
would require $m(2n+1)$ entries.  Both counts are $O(mn)$; the reduction
concerns these arrays and does not imply the same ratio for total memory.

\section{Application to standard two-block ADMM}
\label{sec:admm}

We apply standard two-block ADMM \cite{boyd-admm} to an equivalent
reformulation of the DNN--RLT relaxation~\eqref{eq:mixed-rlt} using auxiliary matrix and cone-image variables. The reformulation
separates the box-cone image from an affine matrix copy containing the normalization and
arrow equations.  We first derive the affine update, then recall the equality-induced
face projection, and finally give the two-block formulation and updates.

\subsection{The affine set and its metric projection}

The affine copy used below need not be constrained to be symmetric because
agreement with the positive semidefinite copy will enforce symmetry.  Its
normalization and arrow equations define
\[
 \mathcal A:=\{W\in\R^{(n+1)\times(n+1)}:W_{00}=1,\quad
 W_{jj}=W_{0j}\ (j=1,\ldots,n)\}.
\]
The remaining subproblem is quadratic minimization over this set.

\begin{lemma}[Projection onto the affine set in the $K$-metric]
\label{lem:arrow-quadratic}
Let $K\in\mathbb S_{++}^{n+1}$ and $B\in\R^{(n+1)\times(n+1)}$.  Let
$e_0,\ldots,e_n$ be the coordinate vectors of $\R^{n+1}$, denote the $j$th
column of $W$ by $w_j$, and denote the $j$th column of $K^{-1}B$ by $z_j$.
Define the constraint coefficient vectors by
\[
 d_0=e_0,\qquad d_j=e_j-e_0\quad(j=1,\ldots,n),
\]
and their right-hand sides by
\[
 \delta_0=1,\qquad \delta_j=0\quad(j=1,\ldots,n).
\]
Then the unique minimizer of
\[
 \min_{W\in\mathcal A}\quad
 \frac12\inner{W}{KW}-\inner{B}{W}
\]
is obtained independently by columns as
\begin{equation}
\label{eq:arrow-rank-one}
 w_j=z_j-K^{-1}d_j
 \frac{d_j^\top z_j-\delta_j}{d_j^\top K^{-1}d_j},
 \qquad j=0,\ldots,n.
\end{equation}
\end{lemma}

\begin{proof}
The affine set $\mathcal A$ imposes the single equation
$d_j^\top w_j=\delta_j$ on column $j$.  Completing the square separates
the problem into
\[
 \min_{d_j^\top w_j=\delta_j}\quad
 \frac12(w_j-z_j)^\top K(w_j-z_j).
\]
Thus $w_j$ is the projection of $z_j$ onto this hyperplane in the norm
$\|v\|_K^2=v^\top Kv$.  The Karush--Kuhn--Tucker equations are
$K(w_j-z_j)+\lambda d_j=0$ and $d_j^\top w_j=\delta_j$.
Solving them gives~\eqref{eq:arrow-rank-one}; its denominator is positive
because $K\succ0$.
\end{proof}

Consequently, this projection first solves $KZ=B$, where the columns
of $B$ are the right-hand sides. Each column $z_j$ of $Z$ is then
corrected along $K^{-1}d_j$ to satisfy $d_j^\top w_j=\delta_j$.  This is
the affine-block update used in the ADMM formulation below.

\subsection{Standard treatment of equality-generated RLT}
\label{sec:equality-face}

We handle equality-generated RLT by the standard one-step facial-reduction
construction~\cite{tuncel-slater,oliveira-wolkowicz-xu,hu-sotirov-wolkowicz,hou-tang-toh}.
Let $\rho=\operatorname{rank}F$, and let
$V\in\R^{(n+1)\times(n+1-\rho)}$ have orthonormal columns spanning
$\ker F$.  The positive semidefinite matrices satisfying the equality-generated
RLT equations $FY=0$ in~\eqref{eq:equality-rlt} form the exposed face
\[
 \mathcal F:=\{Y\succeq0:FY=0\}
 =V\mathbb S_+^{n+1-\rho}V^\top.
\]
Indeed, $Y\succeq0$ and $FY=0$ imply
$\operatorname{range}Y\subseteq\ker F$, so $Y=VRV^\top$ with
$R=V^\top YV\succeq0$; the converse is immediate.  For any matrix $Z$, write $Z_s=(Z+Z^\top)/2$ and $Z_a=(Z-Z^\top)/2$.
Symmetric and skew-symmetric matrices are orthogonal in the Frobenius
inner product. Since $V^\top V=I$, for any symmetric $R$ we have
\[
 \|Z-VRV^\top\|_F^2
 =\|V^\top Z_sV-R\|_F^2
  +\|Z_s\|_F^2-\|V^\top Z_sV\|_F^2+\|Z_a\|_F^2.
\]
Only the first term depends on $R$. Minimizing it over $R\succeq0$ gives
\begin{equation}
\label{eq:face-projection}
 \proj_{\mathcal F}(Z)
 =V\proj_{\mathbb S_+}\!\left(V^\top\frac{Z+Z^\top}{2}V\right)V^\top.
\end{equation}

This standard reduction does not assert that $\mathcal F$ is the minimal
face after the remaining constraints are imposed.  Only the positive semidefinite copy is reduced; the nonnegative and
box-cone copies remain in the original coordinates.
If $p=0$, then $V=I$ and the face projection is the ordinary positive
semidefinite projection.

\subsection{Two-block formulation and updates}

The box-cone representation and the equality-induced face allow us to
separate the relaxation into two blocks with explicit updates.  We keep
the arrow equations in the first block, $W\in\mathcal A$, and introduce
$N$ for the nonnegative matrix copy, $R$ for the reduced positive
semidefinite variable, and $S$ for the box-cone image.  Their cone
constraints define the product cone
\[
 \mathcal K=\{(N,R,S):N\ge0,\ R\succeq0,\
 S\in\mathcal S_m\}.
\]
Taking $(N,R,S)$ as the second block, we write
relaxation~\eqref{eq:mixed-rlt} equivalently as
\begin{equation}
\label{eq:star-consensus}
\begin{aligned}
 \min_{W,(N,R,S)}\quad&\inner{C}{W}\\
 \text{s.t.}\quad&W\in\mathcal A,\qquad (N,R,S)\in\mathcal K,\\
 &\begin{bmatrix}
    W-N\\
    W-VRV^\top\\
    HW-S
   \end{bmatrix}=0.
\end{aligned}
\end{equation}
The stacked equation couples the two blocks linearly: $W=N$ and
$W=VRV^\top$ identify the matrix copies, while $HW=S$ links $W$ to its
box-cone image.  At feasibility, $W=VRV^\top$ makes $W$ symmetric and gives $FW=0$;
hence the columnwise equations defining $\mathcal A$ are the usual arrow
equations.  Conversely, a feasible $Y$ in~\eqref{eq:mixed-rlt} gives a
feasible tuple by taking $W=N=Y$, $R=V^\top YV$, and $S=HY$.  Thus
\eqref{eq:star-consensus} is an equivalent formulation in standard two-block
ADMM form.

We suppress iteration indices and use a superscript $+$ to denote an
updated value; symbols without this superscript denote the current values.
Assign a common penalty $\mu>0$ to all three
coupling equations.  Let $\lambda_N$, $\lambda_P$, and $\lambda_\star$ be the
Lagrange multipliers for $W-N=0$, $W-VRV^\top=0$, and $HW-S=0$,
respectively.  We use the scaled
multipliers
\[
 D_N=\lambda_N/\mu,\qquad
 D_P=\lambda_P/\mu,\qquad
 D_\star=\lambda_\star/\mu.
\]
For any coupling residual $E$ with multiplier $\lambda$ and
$D=\lambda/\mu$, completing the square gives
\[
 \langle\lambda,E\rangle+\frac{\mu}{2}\|E\|_F^2
 =\frac{\mu}{2}\|E+D\|_F^2-\frac{\mu}{2}\|D\|_F^2.
\]
Applying this identity to the three coupling equations and dropping the
terms independent of the primal variables yields the expression below.
Keeping the constraints $W\in\mathcal A$ and
$(N,R,S)\in\mathcal K$ explicit, the scaled augmented
Lagrangian for the three coupling equations, up to terms independent of the
primal variables, is
\[
 \mathcal L_\mu=\inner{C}{W}+\frac{\mu}{2}\bigl(
 \|W-N+D_N\|_F^2+\|W-VRV^\top+D_P\|_F^2
 +\|HW-S+D_\star\|_F^2\bigr).
\]
Each iteration consists of the following three steps.
\begin{enumerate}
\item \textbf{First block update.}

With the second block and scaled multipliers fixed, $W^+$ minimizes
\begin{equation}
\label{eq:admm-affine-quadratic}
 \min_{W\in\mathcal A}\quad
 \frac12\inner{W}{KW}-\inner{B}{W},
\end{equation}
where $\inner{W}{KW}=\operatorname{tr}(W^\top KW)$.  The common coefficient is
\begin{equation}
\label{eq:star-K}
 K=2I+H^\top H
\end{equation}
and the matrix of right-hand sides is
\begin{equation}
\label{eq:star-rhs}
B=N-D_N+VRV^\top-D_P+H^\top(S-D_\star)-C/\mu.
\end{equation}

The image coupling is left multiplication by $H$, whose Frobenius
adjoint is left multiplication by $H^\top$. Its quadratic penalty
therefore contributes $H^\top H$ to $K$, acting independently on each
column. Other scalar formulations may also apply their operators
implicitly, so comparisons with explicit constraint matrices concern a specified implementation choice.
The matrix in~\eqref{eq:star-K} is positive definite, so
Lemma~\ref{lem:arrow-quadratic} applies to~\eqref{eq:admm-affine-quadratic}.
A factorization of $K$ is computed once and reused to obtain all $z_j$ in
the matrix equation $KZ=B$ at each iteration.  The vectors
$K^{-1}d_j$ and denominators $d_j^\top K^{-1}d_j$ are also precomputed.

\item \textbf{Second block update.}

With $W^+$ fixed, minimization of $\mathcal L_\mu$ over $(N,R,S)$
separates into three projections.  Writing
$\operatorname{sym}(Z)=(Z+Z^\top)/2$, the updates are
\[
\begin{aligned}
 N^+&=\max\{W^++D_N,0\},\\
 R^+&=\proj_{\mathbb S_+^{n+1-\rho}}
   \bigl(V^\top\operatorname{sym}(W^++D_P)V\bigr),\\
 S^+&=\proj_{\mathcal S_m}(HW^++D_\star).
\end{aligned}
\]
By~\eqref{eq:face-projection}, the $R$-update gives
$VR^+V^\top=\proj_{\mathcal F}(W^++D_P)$, implementing the face
projection in reduced coordinates.  The box-cone projection is evaluated
row by row using Proposition~\ref{prop:star-projection}.

\item \textbf{Multiplier update.}
Update the scaled multipliers using the coupling residuals:
\[
 D_N^+=D_N+W^+-N^+,\qquad
 D_P^+=D_P+W^+-VR^+V^\top,\qquad
 D_\star^+=D_\star+HW^+-S^+.
\]
\end{enumerate}
\paragraph{Computational cost.}
After forming and factorizing $K$, each iteration costs $O(n^3)$ for
the affine solves, $O(mn^2)$ for multiplication by $H$ and $H^\top$,
$O(n^2)$ for the nonnegative projection, $O(n^3)$ for the positive
semidefinite projection, and $O(mn\log n)$ for the box-cone projections.
Thus both ADMM subproblems use only linear solves and explicit
projections; no constrained QP solver is needed.

\paragraph{Convergence.}
Assuming that \eqref{eq:star-consensus} admits a primal--dual
saddle point, standard two-block ADMM with exact subproblem solutions
and a fixed positive penalty converges. The over-relaxed variant with
parameter $\alpha\in(0,2)$ is covered by
\cite[Theorem~8]{eckstein-bertsekas}.

\paragraph{Valid lower bounds.}
The ADMM implementation separately computes residual-corrected dual bounds using
standard weak-duality postprocessing.  The SDP--RLT comparison tables report returned primal iterate objective
values for both methods, rather than these bounds or SDPNAL+'s returned
dual objectives. These tabulated values should therefore not be
interpreted as certified bounds.

\section{Extensions and related constraint systems}
\label{sec:stronger-rlt}

The shared-bound structure also occurs in related constraint systems.
Full McCormick constraints admit additional box-cone factors, whereas
products between inequality slacks require a different conic image.
We describe these two RLT extensions and their implications for the ADMM
application.

\subsection{Full McCormick constraints}
\label{sec:mccormick-row-star}

The full off-diagonal McCormick system~\cite{mccormick} can be represented
by box cones generated from both bounds $0\le x_i\le1$.  This representation
also enforces entrywise nonnegativity, so the ADMM formulation no longer
needs a separate nonnegative matrix copy.  Define
\[
 H_{\mathrm{box}}=
 \begin{bmatrix}0&I\\\one&-I\end{bmatrix}
 \in\R^{2n\times(n+1)}.
\]
The box bounds are used only as additional RLT generators in the lifted
operator.  Their image is
\begin{equation}
\label{eq:box-row-image}
 H_{\mathrm{box}}Y=
 \begin{bmatrix}x&X\\\one-x&\one x^\top-X\end{bmatrix}.
\end{equation}

\begin{proposition}[McCormick inequalities as box-cone constraints]
\label{prop:mccormick-row-star}
Suppose that $Y\succeq0$ satisfies the normalization and arrow equations.
Then the full off-diagonal McCormick system
\begin{equation}
\label{eq:full-mccormick}
 \max\{0,Y_{0i}+Y_{0j}-1\}\le Y_{ij}
 \le\min\{Y_{0i},Y_{0j}\}
 \quad(1\le i<j\le n)
\end{equation}
is equivalent to
\begin{equation}
\label{eq:mccormick-row-star}
 H_{\mathrm{box}}Y\in\mathcal S_{2n}.
\end{equation}
Either condition implies $Y\ge0$.
\end{proposition}

\begin{proof}
Write $x_i=Y_{0i}$ and $X_{ij}=Y_{ij}$.  The arrow equations and positive
semidefiniteness imply $0\le x_i\le1$.  The box-cone constraints in
\eqref{eq:box-row-image} are therefore
\[
 0\le X_{ij}\le x_i,\qquad
 0\le x_j-X_{ij}\le1-x_i
 \quad(i,j=1,\ldots,n).
\]
For $i\ne j$, these are precisely~\eqref{eq:full-mccormick}; for $i=j$,
they follow from $X_{ii}=x_i$ and $0\le x_i\le1$.
Nonnegativity of $x$ and $X$, together with $Y_{00}=1$, gives $Y\ge0$.
\end{proof}

Accordingly, remove $N$ and the equation $W-N=0$ from
\eqref{eq:star-consensus}, and replace $HW=S$ by
$H_{\mathrm{mc}}W=S_{\mathrm{mc}}$, where
\[
 H_{\mathrm{mc}}=\begin{bmatrix}H\\H_{\mathrm{box}}\end{bmatrix},
 \qquad S_{\mathrm{mc}}\in\mathcal S_{m+2n}.
\]
The second block now contains only the positive semidefinite and box-cone
copies.  With the common penalty $\mu$, the affine subproblem has coefficient
\begin{equation}
\label{eq:mccormick-K}
 K_{\mathrm{mc}}=I+H_{\mathrm{mc}}^\top H_{\mathrm{mc}}
\end{equation}
and matrix of right-hand sides
\begin{equation}
\label{eq:mccormick-rhs}
 B_{\mathrm{mc}}=VRV^\top-D_P
 +H_{\mathrm{mc}}^\top(S_{\mathrm{mc}}-D_{\mathrm{mc}})-C/\mu,
\end{equation}
where $D_{\mathrm{mc}}$ is the scaled multiplier for the augmented box-cone
coupling.  Lemma~\ref{lem:arrow-quadratic} still applies, with one
factorization reused across columns and iterations.  The remaining updates
are the positive semidefinite and box-cone projections and their multiplier
updates from Section~\ref{sec:admm}, using $H_{\mathrm{mc}}$ in place of $H$.

Full McCormick constraints fit naturally into the box-cone formulation:
they add quadratically many scalar inequalities without increasing the
semidefinite matrix order, and the additional box-cone projections cost
$O(n^2\log n)$ per iteration, within a logarithmic factor of the orthant
projection they replace. The comparison with SDPNAL+ uses the DNN--RLT
relaxation~\eqref{eq:mixed-rlt}, without the full McCormick constraints
described here. Their effect on bound quality and solution time is not evaluated in this study.

\subsection{Products between inequality slacks}
\label{sec:inequality-products}

Products between pairs of inequality slacks give another RLT
family.  At a binary rank-one matrix, the slack vector $s=b-Ax$ satisfies
$ss^\top=HYH^\top$, so linearization gives
\begin{equation}
\label{eq:inequality-product-rlt}
 HYH^\top\ge0
\end{equation}
entrywise.  The diagonal inequalities are already implied by
$Y\succeq0$, leaving at most $\binom{m}{2}$ nontrivial off-diagonal
inequalities.  The full image has $m^2$ entries, and storing it together
with its associated multipliers requires $\Theta(m^2)$ storage.  These
products couple different inequality slacks and do not have the shared-slack
structure of the bound-factor products.  They therefore fall outside the
rowwise box-cone representation used here and are excluded from the
relaxation used in the SDPNAL+ comparison.

\section{Numerical experiments}
\label{sec:experiments}

We evaluate the modeling idea in two settings.  The first solves the fixed
SDP--RLT relaxation~\eqref{eq:mixed-rlt} on 257 binary quadratic programs
and compares our box-cone ADMM implementation with SDPNAL+.  The second uses the
box-cone representation of selected triangle inequalities to compute
bounds within MADAM's Max-Cut branch-and-bound framework on 220 graphs.
The first study measures relaxation solution time; the second measures
complete search time and incumbent quality at the time limit.

\paragraph{\textcolor{black}{AI assistance.}}
\textcolor{black}{The authors used AI tools to polish the manuscript's
language and take full responsibility for its contents.}

\subsection{Fixed SDP--RLT relaxations}
\label{sec:fixed-experiments}

The implementation uses standard over-relaxation
\cite{eckstein-bertsekas}, scaling, and a bounded number of penalty
adjustments.  Elementary preprocessing removes redundant inequalities
and identifies forced-zero cone-image coordinates while preserving the
relaxation and the affine update.  Section~S1 of the supplement specifies
these procedures and the stopping criteria.

\subsubsection{Scope and computational protocol}
\label{sec:numerical-protocol}

Both methods solve relaxation~\eqref{eq:mixed-rlt}.  SDPNAL+
\cite{yang-sun-toh,sun-toh-yuan-zhao} represents the inequality-generated RLT constraints as
individual scalar inequalities.

We select inequality-rich instances with $n\le1000$ satisfying
\begin{equation}
\label{eq:numerical-scope-screen}
 m\ge \frac{n+1}{2}.
\end{equation}
These rules are applied before either solver is run.

Each method has a 3600-second iterative-solve limit and a 300-second
setup cap.  We report setup time, solve time, and their sum; Section~S1
of the supplement specifies the timing components.  Both methods use
one CPU thread, 64 GiB of memory, MATLAB R2025b, and the same AMD EPYC
9654 node class.

Both solvers receive the stopping tolerance $10^{-5}$ and terminate according
to their own stopping systems.  We classify a run as solved when the solver
reports convergence; for SDPNAL+, this means termination code zero.
Because their splittings and residual definitions differ, the reported
aggregate residuals are method-specific diagnostics, not a common accuracy
certificate.  Setup limits, solve limits, memory limits, unsuccessful solver
termination, and software errors are reported separately.  Timing comparisons
use total time on instances with convergence reported by both solvers.
This comparison evaluates the two implementations, including their different
algorithms and stopping systems.  Consequently, the timing advantage over
SDPNAL+ should be interpreted under this protocol and cannot be attributed
solely to the box-cone representation.

\subsubsection{Benchmark formulations and data}
\label{sec:benchmark-data}

We summarize the constraint structures, benchmark sources, and dimensions.
Section~S3 of the supplement gives the complete binary formulations,
objectives, variable indexing, and source-to-model conversion conventions.
All objectives and coefficients are taken from the cited collections.

The collections cover several inequality structures relevant to the box-cone
representation: capacity and conflict constraints in quadratic knapsack
\cite{zhou-qkpcg}; graph constraints in DIMACS maximum clique and fixed-color
feasibility \cite{johnson-trick}; covering constraints in OR-Library set
covering \cite{beasley-orlib} and minimum-weight dominating set
\cite{jovanovic-mwds}; three-cycle constraints in bipartite crossing
minimization and single-row facility layout \cite{polip}; triangle
constraints in clique partitioning \cite{sorensen-letchford}; and assignment
and opening constraints in single-allocation hub median
\cite{ernst-krishnamoorthy} and uncapacitated facility location
\cite{beasley-orlib}.  QPLIB QBL contributes general binary quadratic models
with their original linear row bounds \cite{furini-qplib}.

Table~\ref{tab:benchmark-panel} summarizes the benchmark set.  The
minimum-weight dominating-set archive contains many more admissible instances;
to prevent one collection from dominating the experiment, we select source
seed zero for both weight types, every supplied density, and sizes 500, 800,
and 1000.  This gives 30 instances.  For the other collections, we include
all instances covered by the source-data protocols in Supplement~S3 that satisfy the stated selection
rules.

\begin{table}[!ht]
\centering
\small
\caption{Benchmark set.  Here $n$, $p$, and $m$ are the
numbers of binary variables, equality constraints, and inequality
constraints, respectively; the lifted matrix has order $n+1$.}
\label{tab:benchmark-panel}
\begin{tabular}{@{}lrrrr@{}}
\toprule
collection & instances & $n$ & $p$ & $m$\\
\midrule
quadratic knapsack with conflict graphs & 45 & 100--200 & 0 & 199--3,185\\
QPLIB QBL & 17 & 90--595 & 0 & 480--13,090\\
DIMACS maximum clique & 58 & 28--1,000 & 0 & 72--377,247\\
OR-Library set covering & 6 & 192--715 & 0 & 240--4,095\\
fixed-$k$ DIMACS coloring & 12 & 44--957 & 11--95 & 80--10,560\\
bipartite crossing minimization & 37 & 90--552 & 0 & 480--8,096\\
single-row facility layout & 6 & 190--595 & 0 & 2,280--13,090\\
single-allocation $p$-hub median & 11 & 25--625 & 6--26 & 20--600\\
minimum-weight dominating set & 30 & 500--1,000 & 0 & 500--1,000\\
clique partitioning & 31 & 435--990 & 0 & 12,180--42,570\\
uncapacitated facility location & 4 & 816 & 50 & 800\\
\midrule
total & 257 & 25--1,000 & 0--95 & 20--377,247\\
\bottomrule
\end{tabular}
\end{table}

\subsubsection{Results and discussion}
\label{sec:results}

Tables~\ref{tab:comparison-summary} and~\ref{tab:comparison-status}
report all 514 method--instance outcomes for the 257-instance benchmark set,
including resource limits.  A solved case denotes solver-reported convergence.
The supplement provides the
per-instance objectives, residuals, termination statuses, iteration counts,
and setup and solve times underlying the aggregate results.

\begin{table}[!ht]
\centering
\small
\caption{
Results for the fixed 257-instance benchmark set.  Counts A/S denote box-cone
ADMM and SDPNAL+, respectively.  The last three columns use only jointly
solved instances, with total method times $t_A$ and $t_S$ including setup
and solve.  A ratio above one favors ADMM; a dash means no jointly solved
instance.}
\label{tab:comparison-summary}
\begin{tabular}{@{}lrrrr@{}}
\toprule
collection & solved A/S & both solved & faster A/S & median $t_S/t_A$\\
\midrule
bipartite crossing minimization & 37/23 & 23 & 23/0 & 15.90\\
fixed-$k$ DIMACS coloring & 5/5 & 5 & 3/2 & 1.48\\
DIMACS maximum clique & 32/15 & 15 & 15/0 & 7.43\\
OR-Library set covering & 6/2 & 2 & 2/0 & 34.51\\
quadratic knapsack with conflict graphs & 45/28 & 28 & 28/0 & 13.90\\
QPLIB QBL & 17/9 & 9 & 9/0 & 20.81\\
single-allocation $p$-hub median & 9/11 & 9 & 0/9 & 0.42\\
single-row facility layout & 6/1 & 1 & 1/0 & 34.70\\
minimum-weight dominating set & 20/2 & 2 & 2/0 & 12.47\\
clique partitioning & 3/0 & 0 & -- & --\\
uncapacitated facility location & 0/1 & 0 & -- & --\\
\midrule
total & 180/97 & 94 & 83/11 & 13.23\\
\bottomrule
\end{tabular}
\end{table}

\begin{table}[!ht]
\centering
\small
\caption{
Complete status accounting for all 257 instances per method, using solver-reported
termination.}
\label{tab:comparison-status}
\begin{tabular}{@{}lrrr@{}}
\toprule
method & solved & solver timeout & out of memory\\
\midrule
box-cone ADMM & 180 & 77 & 0\\
SDPNAL+ & 97 & 150 & 10\\
\bottomrule
\end{tabular}
\end{table}

On the 94 jointly solved instances, ADMM has the smaller total time on
83 and SDPNAL+ on 11.  The median ratio $t_S/t_A$ is 13.23.  This is a
conditional comparison on the jointly solved subset, not a speedup over the
whole benchmark set.  Performance varies by collection: ADMM is faster on most jointly
solved instances, while SDPNAL+ is faster on all jointly solved
hub-median instances.

Clique partitioning and uncapacitated facility location remain difficult
for both methods.  Some collection-level timing ratios use only one or
two jointly solved instances.

SDPNAL+ reaches the memory limit on ten large maximum-clique models;
ADMM reaches the solve limit on those same models.  All other
unsolved runs reach the solve limit.

\subsection{Max-Cut within MADAM's branch-and-bound framework}
\label{sec:maxcut-experiments}

The second study tests whether box-cone bounds remain useful when cuts
are selected during a branch-and-bound search.  We use MADAM
\cite{hrga-povh-madam} as the host framework and compare its native
cut-strengthened ADMM bounding step with a box-cone replacement.
MADAM supplies the branching rules, primal heuristics, and outer
cut-management routines in both implementations.  Its initial
interior-point solve of the basic SDP relaxation is also retained.

\subsubsection{Bounding implementation and comparison}

The baseline retains MADAM's triangle, pentagonal, and heptagonal cuts.
The box-cone variant uses only triangle cuts, grouping the selected
inequalities into box-cone constraints.  The resulting splitting admits independent box-cone
projections and an affine update with a diagonal coefficient matrix.
Appendix~\ref{sec:maxcut-box} gives the representation and its connection
with bound-factor RLT.  Cuts are generated iteratively using MADAM's
triangle separation routine.  Whenever a triple is selected, the box
variant completes its four triangle inequalities and groups triples
sharing a base pair.  Thus it does not impose all triangle inequalities
at the start of the search.  Scalar cut multipliers recovered from the
cone multipliers are passed back to MADAM for its cut-management steps.

The configured triangle capacities are 50,000 for the baseline and
200,000 for the box variant, allowing four-inequality completion of
selected triples. These are capacity limits, not the numbers of cuts
retained at every node; the larger capacity can also affect cut retention.
The routine uses the same upper-bound correction
as MADAM \cite{hrga-povh-madam}, after recovering nonnegative scalar cut
multipliers from the box-cone multipliers; Section~S5 of the supplement summarizes the bound computation and gives implementation details.

The comparison retains the
search framework but changes the bounding scheme, cut families, and
cut management; the resulting search trees need not coincide.

\subsubsection{Instances and protocol}

We use all 220 Rudy instances with at most 150 vertices distributed in
the official MADAM repository, at commit
\texttt{a051ba3690aa1f549fae8c5ed68348ab02b0a4e0}.\footnote{Source and instances:
\url{https://github.com/HrgaT/MADAM}.}
There are 10, 30, 90, and 90 instances with 60, 80, 100, and 150 vertices,
respectively.  This is a repository-based selection, rather than a
reproduction of the instance selection in the published MADAM experiment.
The panel includes 24 instances used during development, which we also
separate in the analysis below.

Each method is run three times per instance, giving 1320 searches.
Each search has a six-hour wall-time limit, including process startup
and incumbent output.  Both implementations use native C code, the same
compiler and numerical libraries, three MPI processes (one coordinator
and two workers), one BLAS thread per process, and 64 GiB on the same
AMD EPYC 9654 node class as the first study.

A completed search requires branch-and-bound termination reported by the solver and
independent verification of the returned cut using the original integer
edge weights.  Saved incumbent partitions are verified in the same way
for interrupted searches.  These checks validate cut values; optimality
is determined by the solver's numerical bounding and termination rules.
Timing comparisons use only instances completed in all three repeats by
both methods.  For each such instance, let $t_M$ and $t_A$ be the median
whole-search times for MADAM and box-cone ADMM, respectively; we report
the median of $t_M/t_A$ across instances.  All limits are retained in the
completion counts.  For mutually interrupted searches we compare verified
incumbent lower bounds, without treating root bounds as final search-tree
bounds.

The comparison evaluates the complete implementations; it does not isolate the effect of the representation from the changes in cut selection and bounding.

\subsubsection{Results}

Table~\ref{tab:maxcut-comparison} summarizes the full panel.
Per-instance Max-Cut results are provided in Section~S6 of
the supplement.  The box
variant completes all three repeats on 203 instances and MADAM on 202.
Both complete every instance with at most 100 vertices.  Among the 199
jointly completed instances, the box variant is faster on 128 and MADAM
on 71.  The median ratio $t_M/t_A$ is 1.069, corresponding to a median
time reduction of approximately 6.5\%.  Across individual searches, the
box variant completes 612 and reaches the limit on 48; MADAM completes
608 and reaches the limit on 52.  There are no software or scheduler
failures.  Completed objective values agree across methods and repeats
wherever available.

\begin{table}[!ht]
\centering
\small
\caption{Max-Cut comparison on all 220 Rudy instances.  Completion requires
all three repeats.  Counts M/A denote MADAM and box-cone ADMM.
Timing ratios use median whole-search times on jointly completed instances;
values above one favor the box variant.}
\label{tab:maxcut-comparison}
\begin{tabular}{@{}rrrrrr@{}}
\toprule
vertices & instances & completed M/A & both completed & faster M/A & median $t_M/t_A$\\
\midrule
60  & 10 & 10/10 & 10 & 4/6 & 1.019\\
80  & 30 & 30/30 & 30 & 11/19 & 1.040\\
100 & 90 & 90/90 & 90 & 28/62 & 1.123\\
150 & 90 & 72/73 & 69 & 28/41 & 1.055\\
\midrule
all & 220 & 202/203 & 199 & 71/128 & 1.069\\
\bottomrule
\end{tabular}
\end{table}

The advantage varies across graph families.  The median timing ratio is
1.329 on the 17 jointly completed \texttt{pw05} instances, but is slightly
below one on \texttt{g05} and \texttt{w09} (0.980 and 0.985).
On the 196 instances not used during development, both methods complete
all repeats on 180 instances; among their 177 timing pairs, the box
variant is faster on 114, with median ratio 1.067.  The modest timing
advantage therefore also appears outside the development subset.

Both methods reach the limit in every repeat on 13 instances.
Median verified incumbent values agree on ten of these, favor MADAM on
two, and favor the box variant on one.
On the same 199 jointly completed instances, the box variant has a
larger median search-node count on 180 instances, a smaller count on one,
and the same count on 18. Taking the median node count over three repeats
for each method, the median paired box/MADAM node-count ratio is 2.333.
The modest whole-search timing advantage therefore accompanies larger
search trees on most instances. These counts help interpret the complete
implementations, but do not isolate the cost of the bounding routine
from branching, heuristics, cut management, and other search overhead.
These results illustrate the practical value of the box-cone
representation within an established branch-and-bound solver: using only
triangle cuts, the modified bounding routine achieves comparable completion
rates and a modest median timing advantage over MADAM's baseline implementation,
which also uses pentagonal and heptagonal cuts.

\section{Conclusions and scope}
\label{sec:conclusions}

The bound-factor RLT constraints generated by each original inequality
admit an exact representation in the standard box cone.  Their shared
slack gives a compact matrix image and independent projections, preserving
the relaxation and the order of the semidefinite matrix.  The formulation
applies to arbitrary signed linear rows.

Within standard two-block ADMM, the representation yields explicit
cone projections and an affine solve with a reusable factorization.
The experiments illustrate its use for both fixed SDP--RLT relaxations
and dynamically selected Max-Cut triangles within MADAM's search
framework.  Together, these applications support its usefulness as a
building block for first-order methods. Further directions include exploiting the representation in other
first-order splittings and in low-rank methods, including RiNNAL+
\cite{hou-tang-toh-rinnalplus} and the augmented mixing method
\cite{brosch-schwiddessen-wiegele}.

\section*{Funding}
The work of Hao Hu was supported in part by the Air Force Office of
Scientific Research under Award No.~FA9550-23-1-0508.

\section*{Competing interests}
The authors have no competing interests to declare that are relevant
to the content of this article.

{\small
\color{black}
\hypersetup{urlcolor=black,linkcolor=black,citecolor=black}

}

\appendix

\section{Box-cone representation of Max-Cut triangle inequalities}
\label{sec:maxcut-box}

Let an undirected graph have $N$ vertices and edge weights $w_{ij}$.
Max-Cut seeks a partition maximizing the total weight of edges between
its two parts. Represent a cut by signs $z\in\{-1,1\}^N$.  The correlation matrix $Z=zz^\top$ gives cut value
$\frac12\sum_{i<j}w_{ij}(1-Z_{ij})$, where $w_{ij}=0$ for absent edges.
The basic SDP relaxation imposes $Z\succeq0$ and $\diag(Z)=\one$
\cite{GoemansWilliamson1995}.
For distinct vertices $i,j,k$, its four triangle inequalities are
\begin{equation}
\label{eq:maxcut-triangles}
\begin{aligned}
 Z_{ij}+Z_{ik}+Z_{jk}&\ge-1,&
 Z_{ij}-Z_{ik}-Z_{jk}&\ge-1,\\
 -Z_{ij}+Z_{ik}-Z_{jk}&\ge-1,&
 -Z_{ij}-Z_{ik}+Z_{jk}&\ge-1.
\end{aligned}
\end{equation}
These are equivalent to
\[
 |Z_{ik}+Z_{jk}|\le1+Z_{ij},\qquad
 |Z_{ik}-Z_{jk}|\le1-Z_{ij}.
\]
Thus, if $\mathcal T_{ij}$ is a set of selected third vertices for the
base pair $(i,j)$, all corresponding triangles can be represented by
\begin{equation}
\label{eq:maxcut-box}
 \bigl(1+Z_{ij},(Z_{ik}+Z_{jk})_{k\in\mathcal T_{ij}}\bigr)
 \in\mathcal C_{|\mathcal T_{ij}|},\qquad
 \bigl(1-Z_{ij},(Z_{ik}-Z_{jk})_{k\in\mathcal T_{ij}}\bigr)
 \in\mathcal C_{|\mathcal T_{ij}|},
\end{equation}
where
$\mathcal C_d=\{(t,v)\in\R\times\R^d:t\ge0,\ -t\one\le v\le t\one\}$
is the standard box cone over $[-1,1]^d$. Here, $d=|\mathcal T_{ij}|$.
Each factor shares one bound
across its selected third vertices.  For input $c=(c_0,c_1,\ldots,c_d)$, compute $t^\star$ by the
scalar-root calculation and zero-boundary test of
Proposition~\ref{prop:star-projection}, with $|c_j|$ in place of $c_j$
for $j=1,\ldots,d$. Then project each coordinate $c_j$ onto
$[-t^\star,t^\star]$. The same minimization over the coordinates for fixed $t$ gives this formula.

The connection with bound-factor RLT can be seen by fixing vertex $i$ as
an anchor: set $x_j=(1-Z_{ij})/2$, $x_k=(1-Z_{ik})/2$, and let
$X_{jk}=(1-Z_{ij}-Z_{ik}+Z_{jk})/4$ represent $x_jx_k$.
Then~\eqref{eq:maxcut-triangles} is precisely the McCormick system for
$(x_j,x_k,X_{jk})$.  Different triples may use different anchors;
the implementation retains a single correlation matrix and groups
selected triples by a base pair.

Splitting the images in~\eqref{eq:maxcut-box} from the positive
semidefinite copy gives independent box-cone projections.  In symmetric
off-diagonal coordinates the affine update is especially simple: the
identity $(a+b)^2+(a-b)^2=2a^2+2b^2$ cancels the cross terms from the two
images, leaving a diagonal quadratic coefficient for a fixed triangle
set.  Its entries can be reused until that set changes.  This provides
a second application of the shared-bound modeling principle, adapted to
the selected cuts rather than to the complete RLT system in
Section~\ref{sec:model}.

\end{document}


\maketitle

This supplement accompanies the article named above and reports all 257
benchmark instances, with one outcome for the box-cone alternating direction
method of multipliers (ADMM) and one for SDPNAL+ on each instance.
The tables include all outcomes, including resource limits and unsuccessful
termination. The main article gives the formulations, data sources, and
the interpretation of the overall numerical results.
Sections~S1--S4 concern the SDP--RLT comparison.
Section~\ref{sec:supp-maxcut-bound} describes the Max-Cut upper-bound
computation, and Section~\ref{sec:supp-maxcut} reports the 220-instance
comparison with MADAM.

\section{Protocol and interpretation}

Both methods solve the same doubly nonnegative Shor relaxation of a binary
quadratic program, with equality-generated reformulation-linearization
technique (RLT) constraints and the products of every original inequality
slack with each binary variable and its complement. Off-diagonal McCormick
inequalities and products between inequality slacks are excluded.
For each model, $n$ denotes the number of binary variables, $p$ the number
of equality rows, and $m$ the number of inequality rows; the lifted matrix
has order $n+1$. The fixed selection requires $n\le1000$ and
$m\ge(n+1)/2$. The main article describes the selection and benchmark sources.

All executions in the SDP--RLT comparison used MATLAB R2025b Update 3 on AMD EPYC 9654 nodes
of Clemson University's Palmetto cluster, with one CPU thread and 64 GiB
of memory per task. Compiled C++ MEX routines perform box-cone preprocessing
and projection and construct the explicit SDPNAL+ operator. Equality and inequality rows are normalized to unit Euclidean
norm in augmented form before either method is called.

Each method receives a nominal setup allowance of 300 seconds and, if setup
finishes, a solve allowance of 3600 seconds. Total method time is setup plus
solve. Setup includes structural preprocessing, facial-basis construction,
factorization, and initialization for ADMM; for SDPNAL+ it includes explicit
RLT construction, conversion to solver data, and initialization.
MATLAB startup, warm-up, common loading and normalization, independent
verification, and saving results are excluded. Iterative time checks can
allow an overrun by the final iteration; reported times retain the measured
values. Slurm's 75-minute wall-time limit is an outer safety guard.

Both solvers receive tolerance $10^{-5}$. Convergence is determined by the
solver's native stopping criteria; for SDPNAL+, it means termination code
zero. Solver-reported aggregate residuals have method-specific definitions
and are not a common accuracy certificate. Iteration counts are retained for
reference, not for comparing computational work across methods.

The implementation below uses the notation and two-block
formulation of Section~4 of the main article: $W$ is the affine copy,
$N$ the nonnegative copy, $VRV^\top$ the positive semidefinite copy, and
$S$ the box-cone image.  The scaled multipliers are $D_N$, $D_P$, and
$D_\star$; a superscript $+$ denotes an updated value.

\subsection{Over-relaxation}

Over-relaxation is an optional modification intended to improve
numerical progress; it is not needed to define the two-block method.
The choice $\alpha=1$ recovers the ordinary updates in the main article.
We use the standard over-relaxation described by Boyd et al.
(Section~3.4.3 of the ADMM survey cited in the main article;
\url{https://doi.org/10.1561/2200000016}),
with $\alpha=1.7$. This choice can improve practical convergence, but
does not guarantee fewer iterations on every instance.
For $\alpha\in(0,2)$, define the over-relaxed projection arguments

\[
\begin{aligned}
 U_N&=\alpha W^++(1-\alpha)N+D_N,\\
 U_P&=\alpha W^++(1-\alpha)VRV^\top+D_P,\\
 U_\star&=\alpha HW^++(1-\alpha)S+D_\star.
\end{aligned}
\]

The product-space update is then

\[
N^+=\max\{U_N,0\},\qquad
 R^+=\operatorname{proj}_{\mathbb S_+^{n+1-\rho}}
   \bigl(V^\top\operatorname{sym}(U_P)V\bigr),\qquad
 S^+=\operatorname{proj}_{\mathcal S_m}(U_\star).
\]

The scaled-multiplier updates are

\[
\begin{aligned}
 D_N^+&=U_N-N^+,\\
 D_P^+&=U_P-VR^+V^\top,\\
 D_\star^+&=U_\star-S^+.
\end{aligned}
\]

\subsection{Scaling and penalty adaptation}

\label{sec:supp-implementation-scaling}

The following scaling choices are numerical preconditioners; they do not
change the feasible set or the minimizers.  Each nonzero row of $H$ is divided
by a positive scalar.  The implementation assigns the box-cone residual
the penalty $\mu\beta$, where $\beta>0$.  It uses
$D_\star=\lambda_\star/(\mu\beta)$, and replaces $H^\top H$ and
$H^\top(S-D_\star)$ in the affine update by $\beta H^\top H$ and
$\beta H^\top(S-D_\star)$, respectively.  When $C\ne0$, the implementation
also divides $C$ by $\|C\|_F$; reported objective values are restored to the
original scale.

The experiments use $\alpha=1.7$, initial penalty $\mu_0=0.1$, and
$\beta=\min\{2,\max\{0.5,2(n+1)/\max\{1,m\}\}\}$.
Every 25 iterations, let $r_{\mathrm{con}}$ be the maximum of the three
relative consensus residuals defined in the stopping criteria below.
If $r_{\mathrm{con}}>5r_{\mathrm{stat}}$, we double $\mu$; if
$r_{\mathrm{stat}}>5r_{\mathrm{con}}$, we halve it; otherwise we leave it
unchanged.  At most eight penalty changes are allowed.  Each change from
$\mu_{\mathrm{old}}$ to $\mu_{\mathrm{new}}$ multiplies all scaled multipliers
by $\mu_{\mathrm{old}}/\mu_{\mathrm{new}}$, preserving the unscaled multipliers.

Because only finitely many penalty changes are allowed, an indefinitely
continued iteration eventually has a fixed positive penalty. After its
last change, the convergence result stated in Section~4 of the main
article applies under the same saddle-point and exact-subproblem
assumptions. The fixed weight $\beta>0$ can be absorbed by scaling the
third coupling equation by $\sqrt\beta$, as in Boyd et al., Section~3.4.2.

We initialize $W=e_0e_0^\top$, $N=W$, and all scaled multipliers
at zero.  The positive semidefinite copy is initialized by projection onto
$\mathcal F$, and the box-cone image by the projection used after the
preprocessing described below.  The iteration caps are $2\times10^6$ for
ADMM and $2\times10^5$ for SDPNAL+, in addition to the time limits stated
in the main article.

The two solvers use different iteration definitions and perform
different work per iteration, so an equal iteration cap would not impose
an equal computational budget. The requested tolerance and setup/solve
time allowances are shared; the caps are additional termination limits.

\subsection{Elementary preprocessing for structured rows}

Before running ADMM, we use the arrow equations and entrywise nonnegativity
to remove implied RLT inequalities, identify matrix entries forced to zero, and replace
pairs of opposite RLT inequalities by equalities.  These preprocessing steps
preserve the relaxation exactly and require no sign or sparsity assumption
on $A$ and $b$.  They are optional and preserve the two-block ADMM structure,
but modify the box-cone projections by omitting redundant lower bounds and
fixing selected coordinates to zero.  Using the arrow equation in column $j$
gives

\begin{equation}
\label{eq:supp-compiled-q}
\pi_{rj}(Y)=(b_r-a_{rj})Y_{jj}
       -\sum_{i\ne j}a_{ri}Y_{ij}.
\end{equation}

The coefficients in this expression are inspected using exact sign and
coefficient comparisons on the supplied data. In the experiments, these
data have already been normalized before the solver is called.  The three steps are as follows.

\begin{enumerate}

\item
 \textbf{Remove implied inequalities.}
 If $b_r-a_{rj}\ge0$ and $-a_{ri}\ge0$ for every $i\ne j$, then
 $\pi_{rj}(Y)\ge0$ is implied by $Y\ge0$ and may be omitted.

\item
 \textbf{Fix entries to zero.}
 If $b_r-a_{rj}\le0$ and $-a_{ri}\le0$ for every $i\ne j$, the constraint
 $\pi_{rj}(Y)\ge0$ forces $\pi_{rj}(Y)=0$ and every matrix entry with a
 strictly negative coefficient to zero.

\item
 \textbf{Replace opposite inequalities by an equality.}
 If $\pi_{rj}(Y)\ge0$ and $\pi_{r'j'}(Y)\ge0$ are two RLT
 constraints and $\pi_{rj}(Y)=-\pi_{r'j'}(Y)$ holds identically under the
 arrow equations, then both expressions must equal zero.

\end{enumerate}

\Needspace{18\baselineskip}
\subsubsection{Preprocessed row projections}
\label{sec:supp-preprocessed}

For row $r$, let $J_0$ index image coordinates fixed to zero
and let $J_f$ index coordinates whose lower bounds are redundant, with
$J_0\cap J_f=\varnothing$.  The remaining coordinates retain both bounds.
For a projection input $c=(c_0,c_1,\ldots,c_n)$ and a fixed cap $t\ge0$,
the coordinate minimizers are
\[
 u_j(t)=\begin{cases}
 0,&j\in J_0,\\
 \min\{c_j,t\},&j\in J_f,\\
 \min\{\max\{c_j,0\},t\},&j\notin J_0\cup J_f.
 \end{cases}
 \qquad
 \phi_r(t)=t-c_0-\sum_{j\notin J_0}\max\{c_j-t,0\}.
\]
The optimal cap is zero if $\phi_r(0)\ge0$ and otherwise
is the unique positive zero of $\phi_r$.  Fixed coordinates contribute
only a constant to the projection objective, so they do not enter this
equation.  Removing a lower bound changes the coordinate minimizer but
not its contribution to the cap equation.  These formulas follow by the
same scalar minimization as the box-cone projection in the main paper.

The implementation applies these masks to the cone-image
update $S$.  Matrix entries forced to zero by the preprocessing remain
consequences of the model; they are not separately fixed in the
nonnegative-copy update.  The operator $H$, the normalization and arrow
equations, and the affine coefficient $2I+\beta H^\top H$ are unchanged.
Opposite product inequalities are enforced by fixing their image
coordinates to zero.

\subsection{Stopping criteria}

\label{sec:supp-stopping}

We stop when both primal feasibility and stationarity are satisfied to the
prescribed tolerance.  The diagnostics are evaluated after each multiplier
update, using the scaled objective and row operators described in
Subsection~\ref{sec:supp-implementation-scaling}.  Convert the scaled multipliers to
$\lambda_N=\mu_kD_N$, $\lambda_P=\mu_kD_P$, and
$\lambda_\star=\mu_k\beta D_\star$, where $\beta=1$ in the unweighted method.
Let $T_{\mathcal A}$ be the tangent space of $\mathcal A$, consisting of
matrices $Z$ with $Z_{00}=0$ and $Z_{jj}=Z_{0j}$ for $j=1,\ldots,n$.
The relative stationarity residual is

\[
r_{\mathrm{stat}}=
 \frac{\bigl\|\operatorname{proj}_{T_{\mathcal A}}\bigl(C+\lambda_N+\lambda_P
                  +H^\top\lambda_\star\bigr)\bigr\|_F}
 {\max\{1,\|C\|_F,\|\lambda_N+\lambda_P+H^\top\lambda_\star\|_F\}}.
\]

Stationarity is evaluated directly to account for over-relaxation.

For consensus, the Frobenius norm of each difference $W-N$,
$W-VRV^\top$, and $HW-S$ is divided by the maximum of one and the
Frobenius norms of its two terms.  We also
check relative nonsymmetry, $\|W-W^\top\|_F/\max\{1,\|W\|_F\}$, and the
largest absolute equality violation or positive inequality violation among
the normalization and arrow equations, entrywise nonnegativity, $FW=0$,
and all original inequality and bound-factor RLT constraints.  These model
checks use the row data supplied to the solver, before any additional
internal row scaling; in the experiments, the supplied augmented rows
already have unit Euclidean norm.  ADMM reports convergence when every one
of these residuals, including $r_{\mathrm{stat}}$, is at most the requested
tolerance.

Let $P=VRV^\top\succeq0$ be the semidefinite copy. Then
\[
 \operatorname{dist}_F(W,\mathbb S_+^{n+1})\le\|W-P\|_F.
\]
Thus the normalized consensus residual for $W-P$ controls proximity to
the PSD cone at each stopping check. A separate eigenvalue check on the
returned matrix is performed only after termination.  The stopping test measures approximate satisfaction of the
optimality conditions; it does not provide a guaranteed objective gap.

\paragraph{Residual-corrected dual bounds.}
Computing dual bounds by weak-duality postprocessing is standard.
The present SDP--RLT comparison tables report returned primal iterate
objective values for both methods and do not compare corrected dual bounds.

\clearpage
\section{Status accounting and table conventions}

The following labels occur in the detailed tables:
\begin{description}
\item[converged] The solver reports convergence.
\item[solve limit] The iterative solve exhausted its time allowance.
\item[memory limit] The scheduler terminated the task for exceeding its
memory allocation.
\end{description}

Table~\ref{tab:supp-status} accounts for all 514 outcomes: 504 returned
solver results and 10 scheduler-confirmed memory limits.

\begin{table}[ht]
\centering
\caption{Complete outcome counts; each method has 257 tasks.}
\label{tab:supp-status}
\begin{tabular}{lrr}
\toprule
Outcome & ADMM & SDPNAL+\\
\midrule
Converged & 180 & 97\\
Setup limit & 0 & 0\\
Solve limit & 77 & 150\\
Memory limit & 0 & 10\\
\midrule
Total & 257 & 257\\
\bottomrule
\end{tabular}
\end{table}

Tables~\ref{tab:results-1}--\ref{tab:results-11} give the per-instance results.
In each detailed table, ``eq.'' and ``ineq.'' give $p$ and $m$.
Each caption states whether the original problem is a minimization,
maximization, or feasibility problem. ``Obj.'' is the objective value at
the returned iterate, in the original problem's objective sense. 
These values are not certified bounds. Coloring rows have no objective; feasibility of the SDP relaxation
does not establish an integral coloring.
``Solver res.'' denotes the native aggregate residual. Objectives and
residuals are rounded to six and two significant digits, respectively;
times are rounded to seconds, so zero may denote a subsecond time.
Dashes denote unavailable values after a memory limit. Statuses follow
native termination, not rounded residuals.

\clearpage
\section{Benchmark formulations and conversion conventions}
\label{sec:supp-formulations}

This section specifies the binary programs used to construct the relaxations.
All decision variables in the converted models are binary. The instance identifiers in the
per-instance tables identify the retained source models; the main article
provides the bibliographic references and selection rules. Each model is
converted to $\min\{x^\top Qx+q^\top x+\kappa:Ax\le b,\ Gx=h,\
x\in\{0,1\}^n\}$, where $Q$ is symmetric and $\kappa$ is the source objective
constant. A coefficient $a$ on a single off-diagonal monomial $x_ix_j$ is
stored as $Q_{ij}=Q_{ji}=a/2$; a diagonal monomial $a x_i^2$ is stored as
$Q_{ii}=a$. Repeated monomials are accumulated. Maximization objectives,
including their constants, are negated internally and restored for reporting.
The constant does not enter the ADMM updates.

The relaxation and row normalization are specified in Section~S1.

\subsection{Quadratic knapsack with conflicts}
The quadratic knapsack problem with conflict graphs (QKPCG) files supply item profits $a_i$, pair profits $a_{ij}$ for $i<j$,
weights $w_i$, capacity $B$, and an unordered conflict edge set $E_c$.
The formulation is
\[
 \max\ \sum_{i=1}^n a_i x_i+\sum_{1\le i<j\le n}a_{ij}x_ix_j
 \quad\text{subject to}\quad
 \sum_{i=1}^n w_ix_i\le B,\qquad x_i+x_j\le1\quad(ij\in E_c).
\]
Each pair profit is counted once. There are $1+|E_c|$ inequality
constraints and no equalities. We use the 45 medium-sized instances from
groups 1--9 of the published collection.

\subsection{QPLIB QBL}
For each QPLIB file of type QBL (quadratic objective, binary variables,
and linear constraints), let $a_{ij}$ denote its listed quadratic entries with
$i\ge j$, $c_i$ its linear entries, $c_0$ its constant, and $Lx$ its linear
row expressions. The source objective used in the conversion is
\[
 \frac12\sum_i a_{ii}x_i^2+
 \frac12\sum_{i>j}a_{ij}x_ix_j+\sum_i c_ix_i+c_0,
 \qquad \ell_r\le (Lx)_r\le u_r.
\]
Thus a listed diagonal entry contributes $a_{ii}/2$ to $Q_{ii}$ and a
listed off-diagonal entry contributes $a_{ij}/4$ to each symmetric entry
of $Q$, before any objective sign reversal. A row with finite equal bounds becomes an equality. Otherwise each finite
upper bound gives $(Lx)_r\le u_r$ and each finite lower bound gives
$-(Lx)_r\le-\ell_r$. Infinite sides are omitted. The QBL type specifies
binary variables, a quadratic objective, and linear constraints; no
additional application-specific constraints are introduced.

\subsection{Maximum clique and fixed-color feasibility}
Let a DIMACS graph have vertices $1,\ldots,v$ and unordered edge set $E$.
The maximum-clique model is
\[
 \max\sum_{i=1}^v x_i
 \quad\text{subject to}\quad
 x_i+x_j\le1\quad(1\le i<j\le v,\ ij\notin E).
\]
The DIMACS instances are interpreted as undirected simple graphs.
Exactly one inequality is included per nonedge, so $n=v$ and
$m=\binom v2-|E|$.
Selection uses this nonedge count.

For fixed-color feasibility, $x_{it}$ indicates that vertex
$i\in\{1,\ldots,v\}$ receives color $t\in\{1,\ldots,k\}$. Here $k$
is the color count supplied with the DIMACS instance and included in its
name; the counts used below are not assumed to be optimal. The model asks
whether the following constraints admit a binary solution:
\[
 \sum_{t=1}^k x_{it}=1\quad(i=1,\ldots,v),\qquad
 x_{it}+x_{jt}\le1\quad(ij\in E,\ t=1,\ldots,k).
\]
Variables are ordered vertex first, then color. Here $n=vk$, $p=v$, and
$m=k|E|$. The 12 retained models
use \texttt{david}, \texttt{huck} ($k=11$), \texttt{jean} ($k=10$),
\texttt{myciel3}--\texttt{myciel6} ($k=4,5,6,7$), and
\texttt{queen5\_5}--\texttt{queen9\_9} ($k=5,7,7,9,10$).
No symmetry-breaking constraints are added.
A source-graph coloring witness is checked for each model.

\subsection{Set covering and minimum-weight dominating set}

In the set-covering problem, let $a_{ij}=1$ when source row $i$ lists column $j$, and
zero otherwise, and let $c_j$ be the supplied column cost. The formulation is
\[
 \min\sum_{j=1}^n c_jx_j
 \quad\text{subject to}\quad
 \sum_{j=1}^n a_{ij}x_j\ge1\quad(i=1,\ldots,m).
\]
The OR-Library costs and row lists are used directly; conversion negates
each row to obtain a less-than inequality. The retained instances are
\texttt{scpclr10}--\texttt{scpclr13}, \texttt{scpcyc6}, and \texttt{scpcyc7}.

The minimum-weight dominating-set problem has the same objective and
constraints as set covering, with
$m=n$, $c_j$ equal to the supplied vertex weight, and $a_{ij}=1$ precisely
when $i=j$ or vertices $i$ and $j$ are adjacent. The source connection matrix
already includes its unit diagonal. Both supplied weight types are retained
at source seed zero, every supplied density, and vertex counts 500, 800,
and 1000, giving 30 models. Neither covering family has equality rows.

\subsection{Bipartite crossing minimization and single-row facility layout}
For an ordering of $v$ objects, use variables $x_{ij}$ for $1\le i<j\le v$,
with $x_{ij}=1$ indicating that $i$ precedes $j$. The complete constraint
system is
\[
 0\le x_{ij}+x_{jk}-x_{ik}\le1\qquad(1\le i<j<k\le v).
\]
Each triple contributes two inequality rows. Single-row facility layout
uses one ordering, with $n=\binom v2$ and $m=2\binom v3$.
Bipartite crossing minimization uses separate orderings of the two layers,
of sizes $v_1,v_2$, and applies this system separately in each layer. Thus
$n=\binom{v_1}2+\binom{v_2}2$ and
$m=2\binom{v_1}3+2\binom{v_2}3$. Neither model has equality rows.

The selected ordering models from POLIP, the library for polynomially
constrained mixed-integer programming, minimize a free continuous variable $z$
with one coupling row.
Let $x=(x_1,\ldots,x_n)$ list the binary ordering variables in the order
supplied by the source file. Denote the linear, quadratic, and constant
coefficients in its coupling row by $d_a,d_{ab},d_0$, respectively.

The source coupling row is
\[
 z+\sum_a d_a x_a+\sum_{a\le b}d_{ab}x_ax_b+d_0\ge0.
\]
Let $p(x)$ denote the quadratic polynomial following $z$ in this row.
The constraint is therefore $z+p(x)\ge0$, or equivalently $z\ge-p(x)$.
Here $z$ has no separate bounds and appears in no other constraint.
For any fixed feasible ordering $x$, minimizing $z$ therefore sets
$z=-p(x)$. We eliminate $z$ and minimize $-p(x)$, retaining all ordering
constraints. Thus, the equivalent quadratic objective in the binary ordering variables is
$-\sum_a d_ax_a-\sum_{a\le b}d_{ab}x_ax_b-d_0$.
Here $x_a$ denotes the $a$th component of this vector; $a,b$ index those components rather than
the objects being ordered.
The conversion preserves the supplied variable order and all cycle
constraints, and eliminates only $z$.
The quadratic coefficients are symmetrized by the convention at the start
of this section and the constant $-d_0$ is restored in reported objectives.
Source
models containing nonbinary ordering variables are excluded.

For example, POLIP's
\href{https://polip.zib.de/quadratic_linear_ordering/bip_cross_min.16.16.25.1.php}{\texttt{bip\_cross\_min.16.16.25.1}}
has 240 binary ordering variables and the continuous variable $z$;
eliminating $z$ yields the 240-variable model used here.

\subsection{Single-allocation hub median}
Let $v$ be the number of nodes, $f_{ij}$ the supplied directed flows, and
$d_{ij}$ the Euclidean distance between the supplied node coordinates.
Let $a,b,c$ be the collection, transfer, and distribution cost multipliers,
and let $k$ be the prescribed number of hubs. The variable $x_{ih}$ assigns
node $i$ to hub $h$, and $x_{hh}$ opens hub $h$. The formulation minimizes
\[
 \frac1{1000}\left[
 a\sum_{i,h=1}^v\Bigl(\sum_{j=1}^v f_{ij}\Bigr)d_{ih}x_{ih}
 +c\sum_{i,h=1}^v\Bigl(\sum_{j=1}^v f_{ji}\Bigr)d_{hi}x_{ih}
 +b\sum_{i,j,h,l=1}^v f_{ij}d_{hl}x_{ih}x_{jl}\right]
\]
subject to
\[
 \sum_{h=1}^v x_{ih}=1\ (i=1,\ldots,v),\qquad
 \sum_{h=1}^v x_{hh}=k,\qquad
 x_{ih}\le x_{hh}\textcolor{black}{\quad} (i\ne h).
\]
Variables are ordered by assigned node $i$, then hub $h$.
There are $v^2$ variables, $v+1$ equality rows, and $v(v-1)$ inequality
rows; the tautological links with $i=h$ are omitted.
The factor $1/1000$ follows the source-data cost convention. The retained AP instances are \texttt{5.3}, \texttt{10.2}--\texttt{10.5},
\texttt{20.3}, and \texttt{25.3}; the CAB instances are \texttt{10.3.4},
\texttt{15.3.4}, \texttt{20.3.4}, and \texttt{25.3.4}. Their supplied hub
counts and cost multipliers are used without generating new combinations.

\subsection{Clique partitioning}
For $v$ vertices, CP-Lib supplies a weight $w_{ij}$ for each
undirected edge $\{i,j\}$ of the complete graph.
The variable $x_{ij}$ indicates that $i$ and $j$ lie in the same cluster.
The formulation is
\[
 \max\sum_{1\le i<j\le v}w_{ij}x_{ij}
 \quad\text{subject to}\quad
 \begin{cases}
 x_{ij}+x_{ik}-x_{jk}\le1,\\
 x_{ij}+x_{jk}-x_{ik}\le1,\\
 x_{ik}+x_{jk}-x_{ij}\le1
 \end{cases}\quad(1\le i<j<k\le v).
\]
There is no prescribed number of clusters. Pair variables and weights use
the strict upper triangle in row order. The objective is linear in these
variables: internally $Q=0$ and $q_{ij}=-w_{ij}$. There are $\binom v2$
variables, $3\binom v3$ inequalities, and no equality rows.

\subsection{Uncapacitated facility location}
Let $f$ and $c$ denote the numbers of facilities and customers, {\color{black}respectively.} With supplied opening costs
$o_i$ and assignment costs $a_{ji}$, the formulation is
\[
 \min\sum_{i=1}^f o_i y_i+\sum_{j=1}^c\sum_{i=1}^f a_{ji}z_{ij}
 \quad\text{subject to}\quad
 \sum_{i=1}^f z_{ij}=1\ (j=1,\ldots,c),\qquad
 z_{ij}\le y_i\ (i=1,\ldots,f,\ j=1,\ldots,c).
\]
The OR-Library files \texttt{cap71}--\texttt{cap74} are interpreted as
uncapacitated instances: the capacity and demand fields are ignored, and
the supplied assignment costs are not multiplied by demand. The variable
vector lists all opening variables first, followed by assignments in
order of customer index. Thus $n=f(c+1)$, $p=c$, and $m=fc$.
\clearpage
\begin{landscape}
\section{Per-instance results}
\Needspace{12\baselineskip}
\subsection*{Bipartite crossing minimization}
\normalsize
\setlength{\tabcolsep}{1.5pt}
\begin{longtable}{@{}p{.23\linewidth}rrrllrrrrrr@{}}
\caption{Results for Bipartite crossing minimization. Minimization problems.}\label{tab:results-1}\\
\toprule
instance & $n$ & eq. & ineq. & method & status & obj. & solver res. & iter. & setup & solve & total\\
\midrule
\endfirsthead
\toprule
instance & $n$ & eq. & ineq. & method & status & obj. & solver res. & iter. & setup & solve & total\\
\midrule
\endhead
\texttt{bcm\_10\_10\_10\_1} & 90 & 0 & 480 & ADMM & converged & $0.400941$ & $1.0\mathrm{e}{-5}$ & 726 & 0 & 2 & 2\\*
 & & & & SDPNAL+ & converged & $0.401931$ & $7.8\mathrm{e}{-6}$ & 1049 & 1 & 41 & 42\\
\addlinespace[1.5pt]
\texttt{bcm\_10\_10\_20\_1} & 90 & 0 & 480 & ADMM & converged & $8.07616$ & $9.9\mathrm{e}{-6}$ & 18680 & 0 & 57 & 57\\*
 & & & & SDPNAL+ & converged & $8.07877$ & $1.1\mathrm{e}{-4}$ & 2690 & 2 & 123 & 125\\
\addlinespace[1.5pt]
\texttt{bcm\_10\_10\_30\_1} & 90 & 0 & 480 & ADMM & converged & $43.9015$ & $9.9\mathrm{e}{-6}$ & 11992 & 0 & 39 & 39\\*
 & & & & SDPNAL+ & converged & $43.9000$ & $1.4\mathrm{e}{-5}$ & 11129 & 7 & 261 & 267\\
\addlinespace[1.5pt]
\texttt{bcm\_10\_10\_40\_1} & 90 & 0 & 480 & ADMM & converged & $128.837$ & $1.0\mathrm{e}{-5}$ & 13518 & 0 & 35 & 35\\*
 & & & & SDPNAL+ & converged & $128.834$ & $8.7\mathrm{e}{-5}$ & 7300 & 1 & 226 & 227\\
\bottomrule
\end{longtable}
\Needspace{12\baselineskip}
\subsection*{Clique partitioning}
\normalsize
\setlength{\tabcolsep}{1.5pt}
\begin{longtable}{@{}p{.23\linewidth}rrrllrrrrrr@{}}
\caption{Results for Clique partitioning. Maximization problems.}\label{tab:results-2}\\
\toprule
instance & $n$ & eq. & ineq. & method & status & obj. & solver res. & iter. & setup & solve & total\\
\midrule
\endfirsthead
\toprule
instance & $n$ & eq. & ineq. & method & status & obj. & solver res. & iter. & setup & solve & total\\
\midrule
\endhead
\texttt{cpp\_abr\_cars} & 528 & 0 & 16,368 & ADMM & solve limit & $1500.91$ & $7.4\mathrm{e}{-4}$ & 7778 & 2 & 3600 & 3602\\*
 & & & & SDPNAL+ & solve limit & $483.131$ & $8.2\mathrm{e}{-1}$ & 530 & 2 & 3602 & 3604\\
\addlinespace[1.5pt]
\texttt{cpp\_abr\_cetacea} & 630 & 0 & 21,420 & ADMM & solve limit & $975.681$ & $4.0\mathrm{e}{-3}$ & 5161 & 3 & 3600 & 3603\\*
 & & & & SDPNAL+ & solve limit & $-148.015$ & $5.5\mathrm{e}{-1}$ & 337 & 2 & 3617 & 3619\\
\addlinespace[1.5pt]
\texttt{cpp\_abr\_lung\_cancer} & 496 & 0 & 14,880 & ADMM & solve limit & $3468.83$ & $1.6\mathrm{e}{-4}$ & 8770 & 2 & 3600 & 3602\\*
 & & & & SDPNAL+ & solve limit & $899.480$ & $3.4\mathrm{e}{-1}$ & 431 & 2 & 3602 & 3604\\
\addlinespace[1.5pt]
\texttt{cpp\_abr\_micro} & 780 & 0 & 29,640 & ADMM & solve limit & $960.404$ & $1.1\mathrm{e}{-3}$ & 3156 & 5 & 3601 & 3606\\*
 & & & & SDPNAL+ & solve limit & $5.60376$ & $9.3\mathrm{e}{-1}$ & 149 & 4 & 3636 & 3640\\
\bottomrule
\end{longtable}
\Needspace{12\baselineskip}
\subsection*{DIMACS coloring}
\normalsize
\setlength{\tabcolsep}{1.5pt}
\begin{longtable}{@{}p{.23\linewidth}rrrllrrrrrr@{}}
\caption{Results for DIMACS coloring. Fixed-color feasibility problems; objectives are not applicable.}\label{tab:results-3}\\
\toprule
instance & $n$ & eq. & ineq. & method & status & obj. & solver res. & iter. & setup & solve & total\\
\midrule
\endfirsthead
\toprule
instance & $n$ & eq. & ineq. & method & status & obj. & solver res. & iter. & setup & solve & total\\
\midrule
\endhead
\texttt{coloring\_david\_k11} & 957 & 87 & 4,466 & ADMM & solve limit & -- & $2.7\mathrm{e}{-1}$ & 6391 & 1 & 3600 & 3602\\*
 & & & & SDPNAL+ & solve limit & -- & $6.4\mathrm{e}{-1}$ & 420 & 34 & 3606 & 3640\\
\addlinespace[1.5pt]
\texttt{coloring\_huck\_k11} & 814 & 74 & 3,311 & ADMM & solve limit & -- & $2.4\mathrm{e}{-1}$ & 13447 & 1 & 3600 & 3601\\*
 & & & & SDPNAL+ & solve limit & -- & $1.3\mathrm{e}{0}$ & 557 & 32 & 3611 & 3643\\
\addlinespace[1.5pt]
\texttt{coloring\_jean\_k10} & 800 & 80 & 2,540 & ADMM & solve limit & -- & $2.0\mathrm{e}{-1}$ & 14146 & 1 & 3600 & 3601\\*
 & & & & SDPNAL+ & solve limit & -- & $7.8\mathrm{e}{-2}$ & 564 & 27 & 3617 & 3644\\
\addlinespace[1.5pt]
\texttt{coloring\_myciel3\_k4} & 44 & 11 & 80 & ADMM & converged & -- & $9.8\mathrm{e}{-6}$ & 3483 & 0 & 3 & 3\\*
 & & & & SDPNAL+ & converged & -- & $9.6\mathrm{e}{-6}$ & 564 & 1 & 3 & 4\\
\bottomrule
\end{longtable}
\Needspace{12\baselineskip}
\subsection*{DIMACS maximum clique}
\normalsize
\setlength{\tabcolsep}{1.5pt}
\begin{longtable}{@{}p{.23\linewidth}rrrllrrrrrr@{}}
\caption{Results for DIMACS maximum clique. Maximization problems.}\label{tab:results-4}\\
\toprule
instance & $n$ & eq. & ineq. & method & status & obj. & solver res. & iter. & setup & solve & total\\
\midrule
\endfirsthead
\toprule
instance & $n$ & eq. & ineq. & method & status & obj. & solver res. & iter. & setup & solve & total\\
\midrule
\endhead
\texttt{clique\_MANN\_a27} & 378 & 0 & 702 & ADMM & converged & $131.050$ & $9.9\mathrm{e}{-6}$ & 19179 & 0 & 708 & 708\\*
 & & & & SDPNAL+ & converged & $131.117$ & $2.6\mathrm{e}{-5}$ & 14741 & 2 & 3354 & 3356\\
\addlinespace[1.5pt]
\texttt{clique\_MANN\_a9} & 45 & 0 & 72 & ADMM & converged & $17.0882$ & $9.8\mathrm{e}{-6}$ & 1177 & 0 & 1 & 1\\*
 & & & & SDPNAL+ & converged & $17.0899$ & $9.9\mathrm{e}{-6}$ & 1388 & 1 & 7 & 8\\
\addlinespace[1.5pt]
\texttt{clique\_brock200\_1} & 200 & 0 & 5,066 & ADMM & converged & $27.1987$ & $9.8\mathrm{e}{-6}$ & 2215 & 0 & 109 & 109\\*
 & & & & SDPNAL+ & converged & $27.1958$ & $2.6\mathrm{e}{-5}$ & 3021 & 2 & 2568 & 2569\\
\addlinespace[1.5pt]
\texttt{clique\_brock200\_2} & 200 & 0 & 10,024 & ADMM & converged & $14.1321$ & $1.0\mathrm{e}{-5}$ & 1875 & 1 & 151 & 152\\*
 & & & & SDPNAL+ & solve limit & $14.2501$ & $1.8\mathrm{e}{-3}$ & 1389 & 2 & 3603 & 3605\\
\bottomrule
\end{longtable}
\Needspace{12\baselineskip}
\subsection*{Minimum-weight dominating set}
\normalsize
\setlength{\tabcolsep}{1.5pt}
\begin{longtable}{@{}p{.23\linewidth}rrrllrrrrrr@{}}
\caption{Results for Minimum-weight dominating set. Minimization problems.}\label{tab:results-5}\\
\toprule
instance & $n$ & eq. & ineq. & method & status & obj. & solver res. & iter. & setup & solve & total\\
\midrule
\endfirsthead
\toprule
instance & $n$ & eq. & ineq. & method & status & obj. & solver res. & iter. & setup & solve & total\\
\midrule
\endhead
\texttt{mwds\_t1\_v1000e10000\_0} & 1,000 & 0 & 1,000 & ADMM & converged & $1583.69$ & $9.9\mathrm{e}{-6}$ & 3831 & 1 & 1456 & 1456\\*
 & & & & SDPNAL+ & solve limit & $1585.65$ & $1.3\mathrm{e}{-2}$ & 747 & 2 & 3604 & 3606\\
\addlinespace[1.5pt]
\texttt{mwds\_t1\_v1000e1000\_0} & 1,000 & 0 & 1,000 & ADMM & solve limit & $10480.5$ & $1.0\mathrm{e}{-4}$ & 11981 & 0 & 3600 & 3600\\*
 & & & & SDPNAL+ & solve limit & $12542.1$ & $4.0\mathrm{e}{-2}$ & 1455 & 2 & 3601 & 3603\\
\addlinespace[1.5pt]
\texttt{mwds\_t1\_v1000e15000\_0} & 1,000 & 0 & 1,000 & ADMM & converged & $1027.45$ & $9.6\mathrm{e}{-6}$ & 4355 & 1 & 3010 & 3011\\*
 & & & & SDPNAL+ & solve limit & $1024.92$ & $2.1\mathrm{e}{-2}$ & 527 & 2 & 3606 & 3609\\
\addlinespace[1.5pt]
\texttt{mwds\_t1\_v1000e20000\_0} & 1,000 & 0 & 1,000 & ADMM & converged & $770.205$ & $1.0\mathrm{e}{-5}$ & 4485 & 1 & 2878 & 2879\\*
 & & & & SDPNAL+ & solve limit & $767.671$ & $2.3\mathrm{e}{-2}$ & 451 & 2 & 3607 & 3610\\
\bottomrule
\end{longtable}
\Needspace{12\baselineskip}
\subsection*{OR-Library set covering}
\normalsize
\setlength{\tabcolsep}{1.5pt}
\begin{longtable}{@{}p{.23\linewidth}rrrllrrrrrr@{}}
\caption{Results for OR-Library set covering. Minimization problems.}\label{tab:results-6}\\
\toprule
instance & $n$ & eq. & ineq. & method & status & obj. & solver res. & iter. & setup & solve & total\\
\midrule
\endfirsthead
\toprule
instance & $n$ & eq. & ineq. & method & status & obj. & solver res. & iter. & setup & solve & total\\
\midrule
\endhead
\texttt{set\_covering\_scpclr10} & 210 & 0 & 511 & ADMM & converged & $21.5387$ & $9.8\mathrm{e}{-6}$ & 6629 & 0 & 93 & 94\\*
 & & & & SDPNAL+ & converged & $21.5385$ & $1.0\mathrm{e}{-5}$ & 15985 & 1 & 3513 & 3514\\
\addlinespace[1.5pt]
\texttt{set\_covering\_scpclr11} & 330 & 0 & 1,023 & ADMM & converged & $16.8283$ & $1.0\mathrm{e}{-5}$ & 4188 & 0 & 220 & 220\\*
 & & & & SDPNAL+ & solve limit & $17.0518$ & $6.2\mathrm{e}{-3}$ & 2953 & 2 & 3600 & 3602\\
\addlinespace[1.5pt]
\texttt{set\_covering\_scpclr12} & 495 & 0 & 2,047 & ADMM & converged & $16.6365$ & $1.0\mathrm{e}{-5}$ & 2952 & 1 & 566 & 567\\*
 & & & & SDPNAL+ & solve limit & $16.6710$ & $4.7\mathrm{e}{-3}$ & 696 & 2 & 3602 & 3604\\
\addlinespace[1.5pt]
\texttt{set\_covering\_scpclr13} & 715 & 0 & 4,095 & ADMM & converged & $14.4014$ & $9.8\mathrm{e}{-6}$ & 3681 & 5 & 2823 & 2828\\*
 & & & & SDPNAL+ & solve limit & $17.0907$ & $5.2\mathrm{e}{-1}$ & 151 & 12 & 3623 & 3634\\
\addlinespace[1.5pt]
\texttt{set\_covering\_scpcyc6} & 192 & 0 & 240 & ADMM & converged & $51.2009$ & $9.8\mathrm{e}{-6}$ & 3358 & 0 & 28 & 28\\*
 & & & & SDPNAL+ & converged & $51.2011$ & $4.0\mathrm{e}{-5}$ & 15435 & 6 & 885 & 892\\
\addlinespace[1.5pt]
\texttt{set\_covering\_scpcyc7} & 448 & 0 & 672 & ADMM & converged & $121.078$ & $9.8\mathrm{e}{-6}$ & 5492 & 0 & 330 & 330\\*
 & & & & SDPNAL+ & solve limit & $121.257$ & $3.0\mathrm{e}{-3}$ & 10368 & 6 & 3600 & 3607\\
\addlinespace[1.5pt]
\bottomrule
\end{longtable}
\Needspace{12\baselineskip}
\subsection*{QKPCG}
\normalsize
\setlength{\tabcolsep}{1.5pt}
\begin{longtable}{@{}p{.23\linewidth}rrrllrrrrrr@{}}
\caption{Results for QKPCG. Maximization problems.}\label{tab:results-7}\\
\toprule
instance & $n$ & eq. & ineq. & method & status & obj. & solver res. & iter. & setup & solve & total\\
\midrule
\endfirsthead
\toprule
instance & $n$ & eq. & ineq. & method & status & obj. & solver res. & iter. & setup & solve & total\\
\midrule
\endhead
\texttt{1qkpcg1} & 100 & 0 & 199 & ADMM & converged & $17203.7$ & $9.7\mathrm{e}{-6}$ & 10126 & 0 & 27 & 27\\*
 & & & & SDPNAL+ & converged & $17203.7$ & $1.7\mathrm{e}{-5}$ & 11124 & 1 & 290 & 292\\
\addlinespace[1.5pt]
\texttt{1qkpcg2} & 100 & 0 & 199 & ADMM & converged & $13560.0$ & $1.0\mathrm{e}{-5}$ & 7067 & 0 & 20 & 20\\*
 & & & & SDPNAL+ & converged & $13560.0$ & $1.0\mathrm{e}{-5}$ & 7638 & 2 & 173 & 175\\
\addlinespace[1.5pt]
\texttt{1qkpcg3} & 100 & 0 & 199 & ADMM & converged & $16360.0$ & $1.0\mathrm{e}{-5}$ & 32108 & 0 & 104 & 104\\*
 & & & & SDPNAL+ & converged & $16360.3$ & $9.9\mathrm{e}{-6}$ & 165894 & 1 & 1957 & 1958\\
\addlinespace[1.5pt]
\texttt{1qkpcg4} & 100 & 0 & 199 & ADMM & converged & $20121.8$ & $1.0\mathrm{e}{-5}$ & 19087 & 0 & 64 & 64\\*
 & & & & SDPNAL+ & converged & $20121.6$ & $1.2\mathrm{e}{-5}$ & 16929 & 2 & 516 & 519\\
\bottomrule
\end{longtable}
\Needspace{12\baselineskip}
\subsection*{QPLIB QBL}
\normalsize
\setlength{\tabcolsep}{1.5pt}
\begin{longtable}{@{}p{.23\linewidth}rrrllrrrrrr@{}}
\caption{Results for QPLIB QBL. Instances 5935, 5944, 5962, 5971, and 5980 are maximization problems; the remaining instances are minimization problems.}\label{tab:results-8}\\
\toprule
instance & $n$ & eq. & ineq. & method & status & obj. & solver res. & iter. & setup & solve & total\\
\midrule
\endfirsthead
\toprule
instance & $n$ & eq. & ineq. & method & status & obj. & solver res. & iter. & setup & solve & total\\
\midrule
\endhead
\texttt{qplib\_2315} & 595 & 0 & 13,090 & ADMM & converged & $-44088.0$ & $1.0\mathrm{e}{-5}$ & 6221 & 3 & 2773 & 2776\\*
 & & & & SDPNAL+ & solve limit & $-18869.1$ & $7.5\mathrm{e}{-1}$ & 435 & 2 & 3606 & 3608\\
\addlinespace[1.5pt]
\texttt{qplib\_2357} & 240 & 0 & 2,240 & ADMM & converged & $-732.169$ & $9.9\mathrm{e}{-6}$ & 3377 & 0 & 88 & 88\\*
 & & & & SDPNAL+ & converged & $-732.014$ & $2.1\mathrm{e}{-4}$ & 2204 & 1 & 868 & 869\\
\addlinespace[1.5pt]
\texttt{qplib\_2359} & 306 & 0 & 3,264 & ADMM & converged & $-717.610$ & $1.0\mathrm{e}{-5}$ & 18978 & 0 & 1014 & 1014\\*
 & & & & SDPNAL+ & solve limit & $-717.605$ & $3.9\mathrm{e}{-4}$ & 5232 & 1 & 3600 & 3601\\
\addlinespace[1.5pt]
\texttt{qplib\_3584} & 528 & 0 & 10,912 & ADMM & converged & $-37749.5$ & $1.0\mathrm{e}{-5}$ & 4958 & 1 & 1397 & 1398\\*
 & & & & SDPNAL+ & solve limit & $-4223.86$ & $5.4\mathrm{e}{-1}$ & 559 & 2 & 3628 & 3630\\
\bottomrule
\end{longtable}
\Needspace{12\baselineskip}
\subsection*{Single-allocation p-hub median}
\normalsize
\setlength{\tabcolsep}{1.5pt}
\begin{longtable}{@{}p{.23\linewidth}rrrllrrrrrr@{}}
\caption{Results for Single-allocation p-hub median. Minimization problems.}\label{tab:results-9}\\
\toprule
instance & $n$ & eq. & ineq. & method & status & obj. & solver res. & iter. & setup & solve & total\\
\midrule
\endfirsthead
\toprule
instance & $n$ & eq. & ineq. & method & status & obj. & solver res. & iter. & setup & solve & total\\
\midrule
\endhead
\texttt{phub\_ap\_10\_2} & 100 & 11 & 90 & ADMM & converged & $167491$ & $1.0\mathrm{e}{-5}$ & 10375 & 0 & 20 & 20\\*
 & & & & SDPNAL+ & converged & $167491$ & $9.7\mathrm{e}{-6}$ & 1221 & 1 & 12 & 13\\
\addlinespace[1.5pt]
\texttt{phub\_ap\_10\_3} & 100 & 11 & 90 & ADMM & converged & $136007$ & $9.9\mathrm{e}{-6}$ & 14036 & 0 & 32 & 32\\*
 & & & & SDPNAL+ & converged & $136008$ & $5.5\mathrm{e}{-6}$ & 1662 & 1 & 20 & 21\\
\addlinespace[1.5pt]
\texttt{phub\_ap\_10\_4} & 100 & 11 & 90 & ADMM & converged & $112395$ & $9.9\mathrm{e}{-6}$ & 15266 & 0 & 26 & 26\\*
 & & & & SDPNAL+ & converged & $112395$ & $4.5\mathrm{e}{-6}$ & 1970 & 1 & 15 & 16\\
\addlinespace[1.5pt]
\texttt{phub\_ap\_10\_5} & 100 & 11 & 90 & ADMM & converged & $91106.3$ & $1.0\mathrm{e}{-5}$ & 33747 & 0 & 57 & 57\\*
 & & & & SDPNAL+ & converged & $91105.3$ & $8.5\mathrm{e}{-6}$ & 2324 & 1 & 17 & 18\\
\bottomrule
\end{longtable}
\Needspace{12\baselineskip}
\subsection*{Single-row facility layout}
\normalsize
\setlength{\tabcolsep}{1.5pt}
\begin{longtable}{@{}p{.23\linewidth}rrrllrrrrrr@{}}
\caption{Results for Single-row facility layout. Minimization problems.}\label{tab:results-10}\\
\toprule
instance & $n$ & eq. & ineq. & method & status & obj. & solver res. & iter. & setup & solve & total\\
\midrule
\endfirsthead
\toprule
instance & $n$ & eq. & ineq. & method & status & obj. & solver res. & iter. & setup & solve & total\\
\midrule
\endhead
\texttt{srflp\_am\_2009\_33dept\_set1} & 528 & 0 & 10,912 & ADMM & converged & $48340.5$ & $1.0\mathrm{e}{-5}$ & 4958 & 1 & 1777 & 1778\\*
 & & & & SDPNAL+ & solve limit & $72491.9$ & $4.3\mathrm{e}{0}$ & 444 & 2 & 3603 & 3605\\
\addlinespace[1.5pt]
\texttt{srflp\_am\_2009\_35dept\_set1} & 595 & 0 & 13,090 & ADMM & converged & $54783.0$ & $1.0\mathrm{e}{-5}$ & 6221 & 2 & 2930 & 2932\\*
 & & & & SDPNAL+ & solve limit & $97806.3$ & $4.2\mathrm{e}{-1}$ & 560 & 2 & 3607 & 3610\\
\addlinespace[1.5pt]
\texttt{srflp\_anva\_2008\_25dept\_set1} & 300 & 0 & 4,600 & ADMM & converged & $3735.84$ & $9.7\mathrm{e}{-6}$ & 3289 & 0 & 220 & 221\\*
 & & & & SDPNAL+ & solve limit & $3735.20$ & $5.4\mathrm{e}{-3}$ & 4100 & 1 & 3601 & 3602\\
\addlinespace[1.5pt]
\texttt{srflp\_anva\_2008\_30dept\_set1} & 435 & 0 & 8,120 & ADMM & converged & $6655.18$ & $1.0\mathrm{e}{-5}$ & 3068 & 1 & 680 & 681\\*
 & & & & SDPNAL+ & solve limit & $7718.69$ & $6.7\mathrm{e}{-1}$ & 1017 & 2 & 3601 & 3603\\
\addlinespace[1.5pt]
\texttt{srflp\_heku\_1991\_20dept\_set1} & 190 & 0 & 2,280 & ADMM & converged & $13146.9$ & $1.0\mathrm{e}{-5}$ & 1381 & 0 & 27 & 28\\*
 & & & & SDPNAL+ & converged & $13146.8$ & $8.8\mathrm{e}{-6}$ & 2178 & 2 & 956 & 958\\
\addlinespace[1.5pt]
\texttt{srflp\_heku\_1991\_30dept\_set1} & 435 & 0 & 8,120 & ADMM & converged & $35963.9$ & $1.0\mathrm{e}{-5}$ & 4073 & 1 & 670 & 671\\*
 & & & & SDPNAL+ & solve limit & $36215.6$ & $7.2\mathrm{e}{-2}$ & 1003 & 2 & 3602 & 3604\\
\addlinespace[1.5pt]
\bottomrule
\end{longtable}
\Needspace{12\baselineskip}
\subsection*{Uncapacitated facility location}
\normalsize
\setlength{\tabcolsep}{1.5pt}
\begin{longtable}{@{}p{.23\linewidth}rrrllrrrrrr@{}}
\caption{Results for Uncapacitated facility location. Minimization problems.}\label{tab:results-11}\\
\toprule
instance & $n$ & eq. & ineq. & method & status & obj. & solver res. & iter. & setup & solve & total\\
\midrule
\endfirsthead
\toprule
instance & $n$ & eq. & ineq. & method & status & obj. & solver res. & iter. & setup & solve & total\\
\midrule
\endhead
\texttt{ufl\_cap71} & 816 & 50 & 800 & ADMM & solve limit & $935560$ & $6.0\mathrm{e}{-4}$ & 11980 & 0 & 3600 & 3600\\*
 & & & & SDPNAL+ & solve limit & $946180$ & $4.4\mathrm{e}{-4}$ & 1614 & 18 & 3600 & 3619\\
\addlinespace[1.5pt]
\texttt{ufl\_cap72} & 816 & 50 & 800 & ADMM & solve limit & $973264$ & $1.3\mathrm{e}{-3}$ & 10256 & 0 & 3600 & 3600\\*
 & & & & SDPNAL+ & solve limit & $978402$ & $4.1\mathrm{e}{-4}$ & 3842 & 10 & 3605 & 3616\\
\addlinespace[1.5pt]
\texttt{ufl\_cap73} & 816 & 50 & 800 & ADMM & solve limit & $1.01182\mathrm{e}{6}$ & $1.3\mathrm{e}{-4}$ & 15992 & 0 & 3600 & 3600\\*
 & & & & SDPNAL+ & converged & $1.01085\mathrm{e}{6}$ & $1.0\mathrm{e}{-5}$ & 4397 & 10 & 3601 & 3611\\
\addlinespace[1.5pt]
\texttt{ufl\_cap74} & 816 & 50 & 800 & ADMM & solve limit & $1.03287\mathrm{e}{6}$ & $7.3\mathrm{e}{-4}$ & 10315 & 0 & 3600 & 3601\\*
 & & & & SDPNAL+ & solve limit & $1.03656\mathrm{e}{6}$ & $2.6\mathrm{e}{-4}$ & 2447 & 21 & 3600 & 3621\\
\addlinespace[1.5pt]
\bottomrule
\end{longtable}
\end{landscape}
\clearpage
\section{\texorpdfstring{Max-Cut implementation and bounds}{Max-Cut implementation and bounds}}
\label{sec:supp-maxcut-bound}

\subsection{Upper-bound computation}

We adapt MADAM's upper-bound postprocessing to the dual variables
produced by our box-cone formulation. We recover nonnegative scalar cut
multipliers and apply the eigenvalue correction of Hrga and Povh
(Section~4.3; \emph{Computational Optimization and Applications},
80:347--375, 2021,
\href{https://doi.org/10.1007/s10589-021-00310-6}{https://doi.org/10.1007/s10589-021-00310-6}).
This gives an upper bound in exact arithmetic without requiring ADMM
convergence. The floating-point implementation adds a small rounding
cushion but does not provide an interval-certified bound. MADAM's outer
routine adds the objective contribution of fixed variables before
applying its pruning test.

\subsection{Implementation details}

The box-cone subproblem uses standard over-relaxed ADMM with parameter
$1.7$, explicit symmetric-box projections, and full semidefinite
projections.  The penalty starts at one at each search
node and is carried between cut rounds.  Every 25 iterations it is doubled
or halved when the relative consensus and stationarity residuals differ
by a factor greater than ten, within $[10^{-4},10^4]$, with at most eight
changes for each fixed cut system.  Scaled multipliers are rescaled to
preserve their unscaled values.  Both variants retain MADAM's outer tolerance schedule, starting
at $5\times10^{-2}$ and decreasing by a factor of $0.6$ to $10^{-3}$,
and its ADMM iteration-limit parameter of 2000.

Both codes are separate copies of the official MADAM source.  Shared
changes correct the datatype of a transmitted upper bound, make parameter
parsing robust to tabs, and save native output and incumbent partitions.
The baseline retains the original capacity of 50,000 triangle inequalities;
the box variant permits 200,000 to accommodate four-inequality completion.
This can also affect cut retention.

The shared build uses
GCC 12.3.0, OpenMPI 4.1.6, and OpenBLAS 0.3.25.
\clearpage
\begingroup
\section{Max-Cut per-instance results}
\label{sec:supp-maxcut}
We report all 220 Rudy graphs with at most 150 vertices distributed
in the official MADAM repository; $N$ denotes the number of vertices. M denotes MADAM and A denotes box-cone ADMM. Each method has three repeats with a six-hour
whole-search limit per repeat. Completed counts are out of three.
Median time is shown only when all three repeats complete, rounded
to one decimal place; ``limit'' means none completes and ``mixed''
means one or two complete. Best cut is the largest independently
verified cut value across all three repeats, including interrupted
searches.

For the 199 instances completed in all repeats by both methods,
search-node comparisons use the median node count over three repeats
for each method. The median paired box/MADAM node-count ratio is 2.333;
the box variant uses more nodes on 180 instances, fewer on one, and the
same number on 18. Interrupted-search node totals were not recorded and
are excluded from this comparison. Whole-search times and node counts
do not isolate the cost of the bounding routine.
For the 13 instances on which both methods reach the limit in every
repeat, median verified incumbent values agree on ten instances.  The unequal pairs are 17,852
versus 17,832 on \texttt{pw05\_150.0}, 715 versus 713 on
\texttt{pm1d\_150.6}, and 3141 versus 3142 on \texttt{g05\_150.8},
with MADAM listed first.  These comparisons use the median across repeats; the table
reports the largest cut value across repeats.

\normalsize
\begin{longtable}{@{}lrcrr@{}}
\caption{Per-instance Max-Cut results. M: MADAM; A: box-cone ADMM.}\label{tab:supp-maxcut}\\
\toprule
Instance & $N$ & Completed M/A & Median seconds M/A & Best cut M/A\\
\midrule
\endfirsthead
\multicolumn{5}{l}{\tablename\ \thetable\ (continued)}\\
\toprule
Instance & $N$ & Completed M/A & Median seconds M/A & Best cut M/A\\
\midrule
\endhead
\midrule
\multicolumn{5}{r}{Continued on next page}\\
\endfoot
\bottomrule
\endlastfoot
\texttt{g05\_60.0} & 60 & 3/3 & 1.5/1.5 & 536/536\\
\texttt{g05\_60.1} & 60 & 3/3 & 1.2/1.4 & 532/532\\
\texttt{g05\_60.2} & 60 & 3/3 & 1.9/1.8 & 529/529\\
\texttt{g05\_60.3} & 60 & 3/3 & 1.5/1.1 & 538/538\\
\texttt{g05\_60.4} & 60 & 3/3 & 2.1/2.1 & 527/527\\
\texttt{g05\_60.5} & 60 & 3/3 & 1.3/1.3 & 533/533\\
\texttt{g05\_60.6} & 60 & 3/3 & 1.8/1.6 & 531/531\\
\texttt{g05\_60.7} & 60 & 3/3 & 1.9/1.9 & 535/535\\
\texttt{g05\_60.8} & 60 & 3/3 & 1.9/1.6 & 530/530\\
\texttt{g05\_60.9} & 60 & 3/3 & 2.2/2.7 & 533/533\\
\texttt{g05\_80.0} & 80 & 3/3 & 4.4/4.8 & 929/929\\
\texttt{g05\_80.1} & 80 & 3/3 & 1.8/2.3 & 941/941\\
\texttt{g05\_80.2} & 80 & 3/3 & 2.7/2.4 & 934/934\\
\texttt{g05\_80.3} & 80 & 3/3 & 11.8/11.7 & 923/923\\
\texttt{g05\_80.4} & 80 & 3/3 & 4.8/8.3 & 932/932\\
\texttt{g05\_80.5} & 80 & 3/3 & 3.4/3.6 & 926/926\\
\texttt{g05\_80.6} & 80 & 3/3 & 3.0/3.0 & 929/929\\
\texttt{g05\_80.7} & 80 & 3/3 & 2.5/2.0 & 929/929\\
\texttt{g05\_80.8} & 80 & 3/3 & 4.0/4.2 & 925/925\\
\texttt{g05\_80.9} & 80 & 3/3 & 4.8/5.6 & 923/923\\
\texttt{pm1d\_80.0} & 80 & 3/3 & 3.4/3.1 & 227/227\\
\texttt{pm1d\_80.1} & 80 & 3/3 & 7.4/9.5 & 245/245\\
\texttt{pm1d\_80.2} & 80 & 3/3 & 6.2/4.9 & 284/284\\
\texttt{pm1d\_80.3} & 80 & 3/3 & 2.9/2.5 & 291/291\\
\texttt{pm1d\_80.4} & 80 & 3/3 & 3.4/3.3 & 251/251\\
\texttt{pm1d\_80.5} & 80 & 3/3 & 10.0/10.0 & 242/242\\
\texttt{pm1d\_80.6} & 80 & 3/3 & 5.7/5.8 & 205/205\\
\texttt{pm1d\_80.7} & 80 & 3/3 & 8.1/7.6 & 249/249\\
\texttt{pm1d\_80.8} & 80 & 3/3 & 4.8/3.2 & 293/293\\
\texttt{pm1d\_80.9} & 80 & 3/3 & 2.7/2.4 & 258/258\\
\texttt{pm1s\_80.0} & 80 & 3/3 & 1.4/1.2 & 79/79\\
\texttt{pm1s\_80.1} & 80 & 3/3 & 1.6/1.4 & 85/85\\
\texttt{pm1s\_80.2} & 80 & 3/3 & 3.1/2.7 & 82/82\\
\texttt{pm1s\_80.3} & 80 & 3/3 & 2.2/2.0 & 81/81\\
\texttt{pm1s\_80.4} & 80 & 3/3 & 1.6/1.4 & 70/70\\
\texttt{pm1s\_80.5} & 80 & 3/3 & 1.7/1.7 & 87/87\\
\texttt{pm1s\_80.6} & 80 & 3/3 & 2.0/2.0 & 73/73\\
\texttt{pm1s\_80.7} & 80 & 3/3 & 1.8/1.5 & 83/83\\
\texttt{pm1s\_80.8} & 80 & 3/3 & 2.0/2.0 & 81/81\\
\texttt{pm1s\_80.9} & 80 & 3/3 & 1.4/1.4 & 70/70\\
\texttt{g05\_100.0} & 100 & 3/3 & 31.5/41.5 & 1430/1430\\
\texttt{g05\_100.1} & 100 & 3/3 & 183.8/210.9 & 1425/1425\\
\texttt{g05\_100.2} & 100 & 3/3 & 21.2/26.2 & 1432/1432\\
\texttt{g05\_100.3} & 100 & 3/3 & 61.8/47.7 & 1424/1424\\
\texttt{g05\_100.4} & 100 & 3/3 & 7.2/10.2 & 1440/1440\\
\texttt{g05\_100.5} & 100 & 3/3 & 12.8/18.2 & 1436/1436\\
\texttt{g05\_100.6} & 100 & 3/3 & 23.2/21.9 & 1434/1434\\
\texttt{g05\_100.7} & 100 & 3/3 & 31.8/49.8 & 1431/1431\\
\texttt{g05\_100.8} & 100 & 3/3 & 23.1/30.2 & 1432/1432\\
\texttt{g05\_100.9} & 100 & 3/3 & 44.6/21.4 & 1430/1430\\
\texttt{pm1d\_100.0} & 100 & 3/3 & 43.7/45.5 & 340/340\\
\texttt{pm1d\_100.1} & 100 & 3/3 & 80.3/76.2 & 324/324\\
\texttt{pm1d\_100.2} & 100 & 3/3 & 51.0/45.2 & 389/389\\
\texttt{pm1d\_100.3} & 100 & 3/3 & 17.5/21.0 & 400/400\\
\texttt{pm1d\_100.4} & 100 & 3/3 & 42.8/42.4 & 363/363\\
\texttt{pm1d\_100.5} & 100 & 3/3 & 17.9/20.6 & 441/441\\
\texttt{pm1d\_100.6} & 100 & 3/3 & 13.2/13.7 & 367/367\\
\texttt{pm1d\_100.7} & 100 & 3/3 & 10.3/9.8 & 361/361\\
\texttt{pm1d\_100.8} & 100 & 3/3 & 7.3/7.0 & 385/385\\
\texttt{pm1d\_100.9} & 100 & 3/3 & 21.2/19.0 & 405/405\\
\texttt{pm1s\_100.0} & 100 & 3/3 & 3.0/3.1 & 127/127\\
\texttt{pm1s\_100.1} & 100 & 3/3 & 7.7/7.4 & 126/126\\
\texttt{pm1s\_100.2} & 100 & 3/3 & 2.3/2.7 & 125/125\\
\texttt{pm1s\_100.3} & 100 & 3/3 & 12.4/8.4 & 111/111\\
\texttt{pm1s\_100.4} & 100 & 3/3 & 6.5/4.7 & 128/128\\
\texttt{pm1s\_100.5} & 100 & 3/3 & 4.2/2.5 & 128/128\\
\texttt{pm1s\_100.6} & 100 & 3/3 & 4.4/5.2 & 122/122\\
\texttt{pm1s\_100.7} & 100 & 3/3 & 1.4/1.3 & 112/112\\
\texttt{pm1s\_100.8} & 100 & 3/3 & 4.1/3.2 & 120/120\\
\texttt{pm1s\_100.9} & 100 & 3/3 & 2.7/2.9 & 127/127\\
\texttt{pw01\_100.0} & 100 & 3/3 & 4.0/2.7 & 2019/2019\\
\texttt{pw01\_100.1} & 100 & 3/3 & 6.1/4.4 & 2060/2060\\
\texttt{pw01\_100.2} & 100 & 3/3 & 5.6/3.1 & 2032/2032\\
\texttt{pw01\_100.3} & 100 & 3/3 & 3.8/5.0 & 2067/2067\\
\texttt{pw01\_100.4} & 100 & 3/3 & 2.1/1.6 & 2039/2039\\
\texttt{pw01\_100.5} & 100 & 3/3 & 1.7/1.9 & 2108/2108\\
\texttt{pw01\_100.6} & 100 & 3/3 & 5.9/2.7 & 2032/2032\\
\texttt{pw01\_100.7} & 100 & 3/3 & 3.9/4.1 & 2074/2074\\
\texttt{pw01\_100.8} & 100 & 3/3 & 1.8/1.8 & 2022/2022\\
\texttt{pw01\_100.9} & 100 & 3/3 & 9.4/7.6 & 2005/2005\\
\texttt{pw05\_100.0} & 100 & 3/3 & 83.5/62.9 & 8190/8190\\
\texttt{pw05\_100.1} & 100 & 3/3 & 27.1/21.7 & 8045/8045\\
\texttt{pw05\_100.2} & 100 & 3/3 & 41.6/29.1 & 8039/8039\\
\texttt{pw05\_100.3} & 100 & 3/3 & 11.0/8.1 & 8139/8139\\
\texttt{pw05\_100.4} & 100 & 3/3 & 43.0/33.0 & 8125/8125\\
\texttt{pw05\_100.5} & 100 & 3/3 & 11.6/9.6 & 8169/8169\\
\texttt{pw05\_100.6} & 100 & 3/3 & 65.0/53.8 & 8217/8217\\
\texttt{pw05\_100.7} & 100 & 3/3 & 24.1/15.4 & 8249/8249\\
\texttt{pw05\_100.8} & 100 & 3/3 & 21.0/8.2 & 8199/8199\\
\texttt{pw05\_100.9} & 100 & 3/3 & 22.0/18.1 & 8099/8099\\
\texttt{pw09\_100.0} & 100 & 3/3 & 34.8/24.8 & 13585/13585\\
\texttt{pw09\_100.1} & 100 & 3/3 & 50.2/45.5 & 13417/13417\\
\texttt{pw09\_100.2} & 100 & 3/3 & 17.0/13.9 & 13461/13461\\
\texttt{pw09\_100.3} & 100 & 3/3 & 13.9/11.1 & 13656/13656\\
\texttt{pw09\_100.4} & 100 & 3/3 & 21.0/18.9 & 13514/13514\\
\texttt{pw09\_100.5} & 100 & 3/3 & 61.7/42.1 & 13574/13574\\
\texttt{pw09\_100.6} & 100 & 3/3 & 26.0/18.2 & 13640/13640\\
\texttt{pw09\_100.7} & 100 & 3/3 & 51.9/43.6 & 13501/13501\\
\texttt{pw09\_100.8} & 100 & 3/3 & 27.9/20.1 & 13593/13593\\
\texttt{pw09\_100.9} & 100 & 3/3 & 19.2/16.2 & 13658/13658\\
\texttt{w01\_100.0} & 100 & 3/3 & 3.5/2.0 & 651/651\\
\texttt{w01\_100.1} & 100 & 3/3 & 4.5/2.1 & 719/719\\
\texttt{w01\_100.2} & 100 & 3/3 & 11.6/4.4 & 676/676\\
\texttt{w01\_100.3} & 100 & 3/3 & 3.7/2.5 & 813/813\\
\texttt{w01\_100.4} & 100 & 3/3 & 1.7/1.5 & 668/668\\
\texttt{w01\_100.5} & 100 & 3/3 & 5.0/2.7 & 643/643\\
\texttt{w01\_100.6} & 100 & 3/3 & 5.6/10.5 & 654/654\\
\texttt{w01\_100.7} & 100 & 3/3 & 7.0/3.2 & 725/725\\
\texttt{w01\_100.8} & 100 & 3/3 & 1.6/2.3 & 721/721\\
\texttt{w01\_100.9} & 100 & 3/3 & 1.9/1.9 & 729/729\\
\texttt{w05\_100.0} & 100 & 3/3 & 42.9/37.3 & 1646/1646\\
\texttt{w05\_100.1} & 100 & 3/3 & 19.6/14.9 & 1606/1606\\
\texttt{w05\_100.2} & 100 & 3/3 & 13.4/14.2 & 1902/1902\\
\texttt{w05\_100.3} & 100 & 3/3 & 33.5/38.0 & 1627/1627\\
\texttt{w05\_100.4} & 100 & 3/3 & 50.2/56.1 & 1546/1546\\
\texttt{w05\_100.5} & 100 & 3/3 & 66.0/47.1 & 1581/1581\\
\texttt{w05\_100.6} & 100 & 3/3 & 14.3/11.6 & 1479/1479\\
\texttt{w05\_100.7} & 100 & 3/3 & 13.0/12.3 & 1987/1987\\
\texttt{w05\_100.8} & 100 & 3/3 & 45.0/44.3 & 1311/1311\\
\texttt{w05\_100.9} & 100 & 3/3 & 10.9/8.8 & 1752/1752\\
\texttt{w09\_100.0} & 100 & 3/3 & 40.5/33.5 & 2121/2121\\
\texttt{w09\_100.1} & 100 & 3/3 & 208.5/220.6 & 2096/2096\\
\texttt{w09\_100.2} & 100 & 3/3 & 58.2/73.7 & 2738/2738\\
\texttt{w09\_100.3} & 100 & 3/3 & 87.5/79.7 & 1990/1990\\
\texttt{w09\_100.4} & 100 & 3/3 & 50.4/53.6 & 2033/2033\\
\texttt{w09\_100.5} & 100 & 3/3 & 5.9/4.6 & 2433/2433\\
\texttt{w09\_100.6} & 100 & 3/3 & 14.7/10.9 & 2220/2220\\
\texttt{w09\_100.7} & 100 & 3/3 & 25.4/23.8 & 2252/2252\\
\texttt{w09\_100.8} & 100 & 3/3 & 24.3/21.9 & 1843/1843\\
\texttt{w09\_100.9} & 100 & 3/3 & 31.6/33.9 & 2043/2043\\
\texttt{g05\_150.0} & 150 & 0/0 & limit/limit & 3148/3148\\
\texttt{g05\_150.1} & 150 & 3/3 & 12215.3/12377.6 & 3149/3149\\
\texttt{g05\_150.2} & 150 & 3/0 & 17831.5/limit & 3144/3141\\
\texttt{g05\_150.3} & 150 & 0/3 & limit/19860.8 & 3146/3146\\
\texttt{g05\_150.4} & 150 & 3/3 & 19536.1/11678.4 & 3148/3148\\
\texttt{g05\_150.5} & 150 & 3/3 & 7897.4/5536.7 & 3153/3153\\
\texttt{g05\_150.6} & 150 & 0/0 & limit/limit & 3146/3146\\
\texttt{g05\_150.7} & 150 & 0/0 & limit/limit & 3144/3143\\
\texttt{g05\_150.8} & 150 & 0/0 & limit/limit & 3141/3142\\
\texttt{g05\_150.9} & 150 & 0/0 & limit/limit & 3145/3145\\
\texttt{pm1d\_150.0} & 150 & 0/3 & limit/17778.4 & 717/717\\
\texttt{pm1d\_150.1} & 150 & 3/3 & 10734.5/9502.6 & 687/687\\
\texttt{pm1d\_150.2} & 150 & 3/3 & 1298.0/1491.4 & 656/656\\
\texttt{pm1d\_150.3} & 150 & 3/3 & 178.6/209.5 & 762/762\\
\texttt{pm1d\_150.4} & 150 & 3/3 & 18909.7/15652.8 & 671/671\\
\texttt{pm1d\_150.5} & 150 & 3/3 & 17111.8/15866.2 & 606/606\\
\texttt{pm1d\_150.6} & 150 & 0/0 & limit/limit & 715/713\\
\texttt{pm1d\_150.7} & 150 & 3/3 & 17265.7/16632.6 & 708/708\\
\texttt{pm1d\_150.8} & 150 & 3/3 & 12478.8/11555.2 & 739/739\\
\texttt{pm1d\_150.9} & 150 & 3/3 & 1615.7/1590.8 & 640/640\\
\texttt{pm1s\_150.0} & 150 & 3/3 & 600.1/998.9 & 207/207\\
\texttt{pm1s\_150.1} & 150 & 3/3 & 42.5/73.9 & 225/225\\
\texttt{pm1s\_150.2} & 150 & 3/3 & 1014.4/2486.4 & 211/211\\
\texttt{pm1s\_150.3} & 150 & 3/3 & 462.2/370.4 & 228/228\\
\texttt{pm1s\_150.4} & 150 & 3/3 & 168.8/146.8 & 195/195\\
\texttt{pm1s\_150.5} & 150 & 3/3 & 28.1/43.2 & 217/217\\
\texttt{pm1s\_150.6} & 150 & 3/3 & 2720.5/2866.7 & 231/231\\
\texttt{pm1s\_150.7} & 150 & 3/3 & 881.6/649.1 & 223/223\\
\texttt{pm1s\_150.8} & 150 & 3/3 & 1052.2/1328.0 & 175/175\\
\texttt{pm1s\_150.9} & 150 & 3/3 & 1267.0/810.6 & 174/174\\
\texttt{pw01\_150.0} & 150 & 3/3 & 35.9/45.0 & 4402/4402\\
\texttt{pw01\_150.1} & 150 & 3/3 & 226.0/255.2 & 4376/4376\\
\texttt{pw01\_150.2} & 150 & 3/3 & 112.0/91.1 & 4272/4272\\
\texttt{pw01\_150.3} & 150 & 3/3 & 28.1/21.7 & 4285/4285\\
\texttt{pw01\_150.4} & 150 & 3/3 & 295.2/302.3 & 4297/4297\\
\texttt{pw01\_150.5} & 150 & 3/3 & 226.0/34.1 & 4489/4489\\
\texttt{pw01\_150.6} & 150 & 3/3 & 259.7/161.6 & 4446/4446\\
\texttt{pw01\_150.7} & 150 & 3/3 & 142.2/134.7 & 4425/4425\\
\texttt{pw01\_150.8} & 150 & 3/3 & 79.7/43.9 & 4242/4242\\
\texttt{pw01\_150.9} & 150 & 3/3 & 432.5/779.5 & 4287/4287\\
\texttt{pw05\_150.0} & 150 & 0/0 & limit/limit & 17852/17852\\
\texttt{pw05\_150.1} & 150 & 3/3 & 13281.8/8455.6 & 17780/17780\\
\texttt{pw05\_150.2} & 150 & 3/3 & 4680.4/3521.1 & 17774/17774\\
\texttt{pw05\_150.3} & 150 & 3/3 & 14572.5/11109.1 & 17742/17742\\
\texttt{pw05\_150.4} & 150 & 3/3 & 9354.3/6808.0 & 17735/17735\\
\texttt{pw05\_150.5} & 150 & 0/0 & limit/limit & 17867/17867\\
\texttt{pw05\_150.6} & 150 & 3/3 & 12813.6/7267.4 & 18123/18123\\
\texttt{pw05\_150.7} & 150 & 0/0 & limit/limit & 17828/17828\\
\texttt{pw05\_150.8} & 150 & 3/3 & 618.9/602.4 & 17758/17758\\
\texttt{pw05\_150.9} & 150 & 3/3 & 2824.2/2034.1 & 17735/17735\\
\texttt{pw09\_150.0} & 150 & 3/3 & 3118.2/2550.8 & 30015/30015\\
\texttt{pw09\_150.1} & 150 & 3/3 & 6559.7/7803.3 & 30003/30003\\
\texttt{pw09\_150.2} & 150 & 0/3 & limit/17996.5 & 29914/29914\\
\texttt{pw09\_150.3} & 150 & 3/3 & 1386.0/1353.5 & 30059/30059\\
\texttt{pw09\_150.4} & 150 & 3/3 & 12910.9/12169.8 & 29975/29975\\
\texttt{pw09\_150.5} & 150 & 3/3 & 664.3/878.7 & 30301/30301\\
\texttt{pw09\_150.6} & 150 & 3/3 & 5538.5/5913.7 & 30280/30280\\
\texttt{pw09\_150.7} & 150 & 0/1 & limit/mixed & 30211/30211\\
\texttt{pw09\_150.8} & 150 & 3/2 & 19582.0/mixed & 29954/29954\\
\texttt{pw09\_150.9} & 150 & 3/3 & 8072.0/6493.4 & 29878/29878\\
\texttt{w01\_150.0} & 150 & 3/3 & 408.0/377.9 & 1394/1394\\
\texttt{w01\_150.1} & 150 & 3/3 & 199.6/201.9 & 1422/1422\\
\texttt{w01\_150.2} & 150 & 3/3 & 18.6/33.1 & 1315/1315\\
\texttt{w01\_150.3} & 150 & 3/3 & 111.5/110.6 & 1264/1264\\
\texttt{w01\_150.4} & 150 & 3/3 & 64.6/72.5 & 1194/1194\\
\texttt{w01\_150.5} & 150 & 3/3 & 40.4/26.1 & 1303/1303\\
\texttt{w01\_150.6} & 150 & 3/3 & 264.3/288.1 & 1283/1283\\
\texttt{w01\_150.7} & 150 & 3/3 & 205.6/233.5 & 1260/1260\\
\texttt{w01\_150.8} & 150 & 3/3 & 141.5/81.7 & 1167/1167\\
\texttt{w01\_150.9} & 150 & 3/3 & 30.2/28.4 & 1217/1217\\
\texttt{w05\_150.0} & 150 & 3/3 & 14955.7/13564.2 & 3008/3008\\
\texttt{w05\_150.1} & 150 & 3/3 & 6018.8/5389.4 & 3032/3032\\
\texttt{w05\_150.2} & 150 & 3/3 & 2450.7/3284.9 & 2900/2900\\
\texttt{w05\_150.3} & 150 & 2/3 & mixed/12217.6 & 2742/2742\\
\texttt{w05\_150.4} & 150 & 3/3 & 164.8/233.1 & 2815/2815\\
\texttt{w05\_150.5} & 150 & 3/3 & 11721.9/10362.6 & 2609/2609\\
\texttt{w05\_150.6} & 150 & 3/3 & 6143.6/5774.8 & 2835/2835\\
\texttt{w05\_150.7} & 150 & 3/0 & 11040.2/limit & 2608/2608\\
\texttt{w05\_150.8} & 150 & 3/3 & 2120.1/2562.4 & 2781/2781\\
\texttt{w05\_150.9} & 150 & 3/3 & 6495.1/5935.4 & 2725/2725\\
\texttt{w09\_150.0} & 150 & 3/3 & 10884.0/10628.9 & 4278/4278\\
\texttt{w09\_150.1} & 150 & 0/0 & limit/limit & 4177/4177\\
\texttt{w09\_150.2} & 150 & 3/3 & 1656.1/1758.0 & 4019/4019\\
\texttt{w09\_150.3} & 150 & 0/0 & limit/limit & 3451/3451\\
\texttt{w09\_150.4} & 150 & 3/3 & 14991.7/16082.2 & 3831/3831\\
\texttt{w09\_150.5} & 150 & 3/3 & 1730.5/2326.9 & 3579/3579\\
\texttt{w09\_150.6} & 150 & 3/3 & 2853.2/3605.6 & 3889/3889\\
\texttt{w09\_150.7} & 150 & 0/0 & limit/limit & 3304/3304\\
\texttt{w09\_150.8} & 150 & 0/0 & limit/limit & 3976/3976\\
\texttt{w09\_150.9} & 150 & 3/3 & 11603.3/8543.3 & 3593/3593\\
\end{longtable}
\endgroup